\documentclass[12pt]{article}
\usepackage{amsmath}
\usepackage{graphicx}
\usepackage{enumerate}
\usepackage{natbib}
\usepackage{url} 
 \usepackage{fancyhdr}
 \usepackage{epstopdf}
 \usepackage{multirow,booktabs}
\usepackage{color}
  \usepackage{amsthm,amssymb}
\usepackage{mathrsfs}
\usepackage{indentfirst}
\usepackage{enumitem}
\usepackage{indentfirst}
\usepackage{longtable}
\usepackage{float}
\usepackage{xcolor}
\usepackage{bigstrut}
\usepackage{colortbl}
\usepackage{array}
 \usepackage{url}
\usepackage{chngcntr}
\usepackage{epsfig}
\newcommand{\blind}{0}

\def\theequation{\arabic{section}.\arabic{equation}}
\newtheorem{theorem}{Theorem}[section]
\newtheorem{corollary}{Corollary}[section]

\newtheorem{rema}{Remark}[section]

\counterwithout{equation}{section}

\begin{document}

\def\spacingset#1{\renewcommand{\baselinestretch}%
{#1}\small\normalsize} \spacingset{1.5}
	\date{}

\if0\blind
{
  \title{\bf The Role of Propensity Score Structure in Asymptotic Efficiency of Estimated Conditional Quantile Treatment
Effect}
  \author{Niwen Zhou\\
    School of Statistics, Beijing Normal University\\
    and \\
    Xu Guo\\
   School of Statistics, Beijing Normal University\\
    and \\
    Lixing Zhu\thanks{
    The authors gratefully acknowledge \textit{two grants from the University Grants Council of Hong Kong (HKBU123017/17 and HKBU123028/18) and a NSFC grant (NSFC11671042)}}\hspace{.2cm}\\
   Department of Mathematics, Hong Kong Baptist University}
  \maketitle
} \fi

\if1\blind
{
  \bigskip
  \bigskip
  \bigskip
  \begin{center}
    {\LARGE\bf The role of propensity score structure in asymptotic efficiency of estimated conditional quantile treatment
Effect}
\end{center}
  \medskip
} \fi

%
%

\bigskip
\begin{abstract}
When a strict subset of covariates are given, we propose conditional
quantile treatment effect to capture the heterogeneity of treatment effects via the quantile sheet that is the function of the given covariates and quantile.  
We focus on deriving the asymptotic normality of probability score-based
estimators under parametric, nonparametric and semiparametric structure. We make a systematic study on the estimation efficiency to check the importance of propensity score structure and the essential differences
from the unconditional counterparts. The derived unique properties can answer: what is the general ranking of these estimators? how does the affiliation of the given covariates to the set of covariates of the propensity score affect the efficiency? how does the convergence rate of the estimated propensity score affect the efficiency? and why would semiparametric estimation be worth of recommendation in practice?  We also give a brief discussion on the extension of the methods to handle large-dimensional scenarios and on the estimation for the asymptotic variances. The simulation studies are conducted to examine the performances of these estimators. A real data example is analyzed for illustration and some new findings are acquired.
\end{abstract}

\noindent%
 \vfill

\spacingset{1.5} 



%
%
%
\section{Introduction}

{
Treatment effect is a vital issue in diverse research and
applied fields. {In the  literature, most of existing studies  focus on average treatment effect (ATE) defined by the population mean of potential outcomes as well as quantile treatment effects ($QTE$) by the population quantile. $QTE$ can capture the heterogeneity of  treatment effect.} For example, reducing class sizes may have positive effect on the academic performance of excellent students but opposite on that of the weaker, or vice versa \citep{koenker2017}.
}\cite{doksum1974} and \cite{d1975} defined the $\tau$th $QTE$ as the difference between quantiles of the
two marginal potential outcome distributions. As \cite{firpo2007} commented,  the $\tau$th quantile in this definition
is not exactly equal to   the $\tau$th quantile of the distribution of the potential outcomes
difference unless  the rank
preservation assumption is satisfied. {Yet it is still reasonable and informative in studying  treatment effects and some references include \cite{koenker2005}, \cite{firpo2007}, \cite{zhang2018}.}
\cite{firpo2007} proposed  the $QTE$ and $QTT$ (quantile
treatment effect on the treated) estimator  by minimizing
the expectation of proper weighting check functions, where
the weight is based on the inverse of propensity score.
 {Further,  conditional $QTE$ ($CQTE$) can provide quantile sheet, that is a function of quantile levels and  given covariates,
to examine one more type of  heterogeneity:  $QTE$ in specific
subpopulations decided by the covariates, which can  reflect their influence at different quantiles.} {Such a heterogeneity is also informative, see, e.g. \cite{wager2018} and \cite{luo2019}.}

\begin{figure}[htbp]
\vspace*{-10pt}
  \centering
  \includegraphics[scale=0.65]{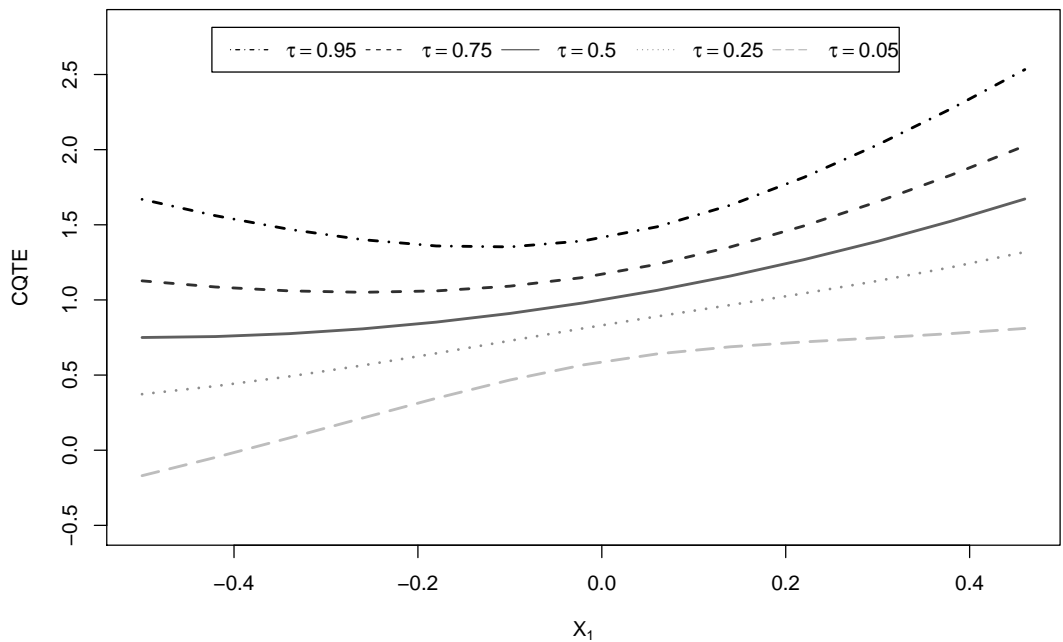}
   \vspace*{-20pt}
 \caption{\scriptsize{The quantile sheet of $CQTE(X_1)$ for Model~1 in simulation with $\tau=\{0.05,0.25,0.5,0.75,0.95\}.$}\label{true_model2}}
\end{figure}

In this article, we propose $CQTE$
in a general situation in which the conditioning continuous
covariates $X_1$ form a strict subset of the covariates
$X$ ($X_1 \subsetneq X$). This is a generalization of $CQTE$ proposed by \cite{imbens2009} who referred to it as the
difference of quantiles for two potential outcome
distributions conditional on the whole $X$ to guarantee the unconfoundedness
assumption. By using conditional quantile treatment
effect ($CQTE(X_1)$, hereafter), we can know how the treatment effect change with $X_1$, which can help not only a detailed programme evaluation, but also the investigation on the importance of  $X_1$. Note that the unconfoundedness assumption on $X_1$
may not hold.  Based on this indirect unconfoundedness assumption the technical skills for theoretical development have to be more sophisticated than those for $QTE$.
Also, the quantile sheet of $CQTE(X_1)$ (denoted as $\Delta_\tau(X_1)$), as a function of both $\tau$ and $X_1$, is more informative than conditional average treatment effect ($CATE(X_1)$) proposed by \cite{abrevaya2015} as shown in Figure~\ref{true_model2}.

It is well known that in unconditional cases, with estimated propensity score,  the estimation efficiency for  $ATE$ and $QTE$ can be enhanced. 
See, \cite{hir2003}. This type of estimator is referred as
the inverse probability weighting-based (IPW, hereafter) estimator.
Based on the nonparametrically estimated propensity score, \cite{firpo2007} proposed an IPW type estimator for $QTE$ that can achieve the semiparametric efficiency bound.  
Based on the potential outcome model, we in this paper will construct the
pointwise $CQTE(X_1)$ estimator via minimizing a properly
weighted sum of check functions with the estimated propensity score. Therefore, we will  first
estimate and then asymptotic behaviours of the estimated $CQTE$ when the propensity score is under parametric, nonparametric  and  semiparametric  dimension reduction structure. See the relevant references such as  \cite{Yao2010},
\cite{abrevaya2015} and \cite{guo2018}. 
{

{As under the regularity conditions designed in this paper, all estimators are asymptotically unbiased, we then discuss asymptotic efficiencies  by the asymptotic variances. According to the research for  $ATE$ and $QTE$,  we consider the efficiency bound and efficient estimation construction for the $CQTE$ function $\Delta_\tau(X_1)$ defined in the next section. As pointed out by \cite{kennedy2017}, if we only assume mild smoothness conditions on $\Delta_\tau(X_1)$, there is no existing theory in the literature to derive the efficiency bound and an efficient estimator for $\Delta_\tau(X_1)$ in the sense ATE or QTE shares. Further,  because $\Delta_\tau(X_1)$ is not pathwise differentiable, any estimator cannot achieve $\sqrt{n}$-consistent. However, we can have some information on the estimation efficiency as follows. For brevity, write $OCQTE$, $PCQTE$, $SCQTE$ and $NCQTE$ as the estimators with true, parametric, semiparametric, and nonparametric estimated propensity score respectively.}  Let $A\preceq B$ mean that the asymptotic variance of estimator $A$ is not greater than that of estimator $B$ and $A\cong B$ stand for that $A$ has the same asymptotic variance function as $B$. We have the following.
\begin{enumerate}
\item In general, the asymptotic efficiency of the four estimators has the ranking: $$ NCQTE \preceq SCQTE \preceq PCQTE  \cong OCQTE. $$
{
\item  { When the estimated propensity score has $X_1$ as a strict subset of its true arguments,} $$ NCQTE \preceq SCQTE \preceq PCQTE  \cong OCQTE. $$ When $X_1$ is not a strict subset of its true arguments, $$ NCQTE \cong SCQTE \cong PCQTE  \cong OCQTE. $$} Subsection~2.5.2 presents some more detail.
\item  If the propensity score function is very smooth and the kernel functions and tuning parameters in the nonparametric estimation are selected delicately, $$ NCQTE \cong SCQTE \cong PCQTE  \cong OCQTE. $$

    Subsection~2.5.1 presents the results. It is worthwhile to point out that the research in this scenario basically serves as a theoretical exploration and provides an insight into the nature of CQTE. For practical use, we may not consider such ways to estimate NCQTE and SCQTE to lose their estimation efficiency. 
\item We recommend $SCQTE$ for practical use as it can very much alleviate the curse of dimensionality which is a very serious problem for $NCQTE$, and is robust against model misspecification, particularly, of parametric propensity score structure.
\end{enumerate}
These newly found phenomena 
show the unique properties of $CQTE$ and demonstrate the essential differences from their unconditional counterparts. As is well known,  estimating propensity scores can always enhance, with smaller asymptotic variances, the estimation efficiencies of their unconditional counterparts. Further, the nonparametrically estimated propensity score can make a better efficiency than the parametric/semiparametric one. Relevant references include Hirano et al. (2003), Guo et al. (2018) and Liu et al. (2018).

{
It should be mentioned that some parts of this research are extensions, but not trivial, of existing works, 
$NCQTE$ is an extension of the procedure of \cite{firpo2007} from $QTE$ to $CQTE$. Since $CQTE$ is a function of the given convariates $X_1$, 
it makes the asymptotic analysis essentially different from that of $QTE$.  $PCQTE$ and $NCQTE$ also extend, with more information, the approach of \cite{abrevaya2015} from $CATE$. But the unsmoothness of the quantile loss function causes the asymptotic analysis  more complex than that for $CATE$. The new $SCQTE$ has a very important feature of dimension reduction nature in estimating propensity score. This feature can simultaneously alleviate the risk of mis-specification and the curse of dimensionality.

}
{
The rest of the paper is organized as follows. In Section~2,
we introduce the estimation procedures for $CQTE(X_1)$ and
investigate their asymptotic properties.   Subsection~2.4 is devoted to give further results about the three estimators, and Subsection~2.5 to present more detailed results about the role of convergence rate of estimated propensity score and the role of the affiliation of $X_1$ to the set of the arguments of propensity score.
Section~3 contains
some numerical studies to examine the performance of the three estimators. In Section~4, we apply
our methods to analyze a real data set for illustration and find some phenomena, which $CATE$ cannot obtain.
Section~5 contains some conclusions and remarks about more general models and the reason why we in this paper do not include the investigation for the other two basic methodologies: potential outcome regression and doubly robust estimation. Due to the space limitation, all the
technical proofs are relegated to the supplementary material.

}

%
%

\section{Estimation procedures and asymptotic properties}
\subsection{Definition and preparation }

Let $D$ be the indicator variable of treatment and $Y$ the
outcome.  $D_i=0,1$ respectively means the $i$th individual
does not receive or receives treatment. Denote
the corresponding potential outcome as $Y_i(0)$ or $Y_i(1)$ and  write the observed
outcome as
$Y_i=D_iY_i(1)+(1-D_i)Y_i(0).$  Let $X$ be a $k$-dimensional
vector of covariates with $k\geq2 $ and $X_1\in R^l$ be a subvector of $X\in R^k$ with $1\leq l<k$. Write $p(X)$ as the propensity score
$E(D\mid X)$.  Further assume that
$(X_i,Y(1)_i,Y(0)_i,D_i), i=1,\cdots,n,$ are independent identically distributed (i.i.d.) random vectors.
Let $\tau$ be a real
value in $(0,1)$ and the $CQTE$ function
$\Delta_\tau(x_{10})=q_{1,\tau}(x_{10})-q_{0,\tau}(x_{10})$ with
\\
  \centerline{$q_{j,\tau}(x_{10})=\inf_a
  E[\rho_{\tau}(Y(j)-a)\mid X_{1}=x_{10}],j=0,1.$}
 Here  $\rho_{\tau}(u)=u(\tau-\mathbb I(u<0))$ is the check
   function, $\mathbb I(\cdot) $ is an  indicator function, and
 {$x_{10}\in \Omega$ with $\Omega$ containing all the interior points of the support of $X_1$.}
  Denote the conditional
 distribution of $Y(1)\mid X_1$ and of $Y(0)\mid X_1$ as
 $F_{1}(Y(1)\mid X_1)$ and $F_{0}(Y(0) \mid X_1)$ respectively.

   To introduce the estimation
procedure and theoretical results smoothly,
 the following assumptions are required.
{
\begin{itemize}
 \item Assumption 1: (Strong ignorability)

 (i) Unconfoundedness: $(Y(0),Y(1))\perp D\mid X$.

 (ii) Common support: For some very small $c>0$,
 $c<p(X)<1-c$.

 \item Assumption 2: (Conditional quantile function) For any $\tau \in(0,1)$, $j=0,1$ and $x_{10}\in \Omega$,

  (i) $q_{j,\tau}(x_{10})$ is the unique $\tau$th
   conditional quantile of $Y(j)\mid X_1=x_{10}$.

 {(ii)  $q_{j,\tau}(X_1)$,$j=0,1$ is $s_1\geq 2$ times continuously differentiable for $X_1$.

 (iii) The conditional distribution function $F_j(y_j\mid X_1)$ and density function $f_j(y\mid X_1)$ are bounded and uniformly continuous in $y_j$ for $X_1$.}
 \end{itemize}

Assumption 1 is commonly used,
see e.g.,  \cite{rosenbaum1983}. 
Part(i) of Assumption 1 implies that
the observed vector $X$ can fully control for any endogeneity in treatment choice and 
part(ii) of Assumption 1 means that there is overlap between
the supports  of the conditional distributions of $X$ given $D=0$ and $D=1$ respectively. Part(i) of Assumption 2
  guarantees the identifiability of $q_{j,\tau}(x_{10})$, $j=0,1.$ It follows that
$q_{1,\tau}(x_{10})=F_1^{-1}(\tau\mid X_1=x_{10})$ and $q_{0,\tau}(x_{10})=F_0^{-1}(\tau\mid X_1=x_{10}),\forall \tau\in(0,1),x_{10}\in \Omega.$
Further, part(ii) of Assumption 2 is required to ensure the function smoothness  which will be used, particularly for nonparametric-based estimations.
}

Under the unconfoundedness assumption, $q_{j,\tau}(x_{10})$
for $j=0, 1$ can be also rewritten as
 \begin{eqnarray}
 \begin{split}
 q_{0,\tau}(x_{10})&=\inf_{a_0}
 E\left[\frac{1-D}{1-p(X)}\rho_{\tau}(Y-a_0)\mid
 X_{1}=x_{10}\right],\\\label{est_2.1}
 q_{1,\tau}(x_{10})&=\inf_{a_1}
 E\left[\frac{D}{p(X)}\rho_{\tau}(Y-a_1)\mid
 X_{1}=x_{10}\right],
 \end{split}
 \end{eqnarray}
 where $p(x)=P(D=1\mid X=x).$

  Note that
estimating $\Delta_\tau(x_{10})= q_{1,\tau}(x_{10})- q_{0,\tau}(x_{10})$ does not involve estimating
the conditional distributions $F_1(Y(1)\mid X_1)$ and $F_0(Y(0)\mid X_1)$ that can be
nonparametric. {Thus, we can estimate $\Delta_\tau(x_{10})$ in a
simpler
 manner.}

 After having the estimation of $p(X)$, we then estimate  $q_{j,\tau}(x_{10}),j=0,1$
separately by a  nonparametric method $\hat q_{j,\tau}(x_{10}),j=0,1$ such as local linear smoother (e.g. \citep{fan1996}) and deriving
asymptotically linear representations of $\hat
q_{j,\tau}(x_{10}),j=0,1$ and then of $\hat {\Delta}_\tau
(x_{10})=\hat{q}_{1,\tau}(x_{10})-\hat{q}_{0,\tau}(x_{10}).$ The following
subsections present the estimations and theoretical results.

\subsection{ Three different estimators of $p(X)$}

 If $p(X)=\pi(X,\beta)$ is known up to some unknown parameters $\beta$ such as the
popular logistic model or probit model, we then need to estimate $\beta$. If we do
not have such a prior information on its structure, a
nonparametric estimation is required such as
the Nadaraya-Watson (N-W) estimation.
 {
 Furthermore, when it has a
semiparametric structure: $p(X)=p(\alpha^{\top}X)$, where both the function $p(\cdot)$  and the $k\times q$ orthonormal matrix $\alpha$ are unknown with $q\leq k$.
From the definition of $p(X)$, we can see that
the information about $D$ from $X$ can be completely captured by the projected
variables $\alpha^{\top}X$. Thus, we can use the following
conditional independence to present the above semiparametric
structure:
\begin{eqnarray}\label{cons_1}
\mbox{Constraint~1}:\ \ D\perp X \mid \alpha^{\top} X.
\end{eqnarray}
It follows that $(Y(0),Y(1))\perp D\mid \alpha^{\top}X$.  Note that (\ref{cons_1}) still holds if we replace $\alpha$ by any
$\alpha C$, where $C\in R^{q\times q}$ is any nonsingular matrix.
In general, the matrix $\alpha$ can only be identifiable up to a rotation matrix $C$. Thus, under this dimension reduction framework, \cite{li2018} pointed out that the identifiable parameter in (\ref{cons_1}) is $\alpha C$ or in the other words, the space spanned by the columns of $\alpha$. In the literature, various methods have been proposed to estimate this space including sliced inverse regression (SIR, \citet{li1991}), and minimum average variance estimation (MAVE, \cite{xia2002}, \cite{xia2007}).
{As for determining the structural dimension $q$, several eigen-decomposition-based  methodologies have been proposed in the literature, such as the sequential test methods \citep{li1991} and the BIC-type methods \citep{zhu2006}. For ease of exposition, we assume the dimension $q$ of $\alpha$ is given.} This semiparametric dimension reduction structure can not only alleviate the curse of dimensionality, but also  maintain the model interpretation and flexibility simultaneously to greatly avoid model mis-specification.

}

%


The three  estimators of $p(X)$ in
parametric,
nonparametric and semiparametric scenarios are respectively as
\begin{eqnarray}\label{semi+equ}
\begin{split}
&\hat p(X_i)=\pi(X_i,\hat \beta),\quad \hat
\beta=\mathop{\arg\max}_{\beta}\sum_{i=1}^n(D_i \log
\pi(X_i,\beta)+(1-D_i)(1-\log \pi(X_i,\beta) ); \\
&\hat p(X_i)= {\frac{1}{n{h_0^k}}\sum_{j:j\neq i}D_iL
\left(\frac{X_j-X_i}{h_0}\right)}\big/{\frac{1}{nh_0^k}
\sum_{j:j\neq i}L\left(\frac{X_j-X_i}{h_0}\right)};   \\
&\hat p(X_i)=\hat p(\hat\alpha^{\top}
X_i)= {\frac{1}{nh_2^q}\sum_{j:j\neq i}D_iH\left(\frac{\hat\alpha
^{\top}X_j-\hat\alpha
^{\top}X_i}{h_2}\right)}\big/{\frac{1}{nh_2^q}\sum_{j:j\neq
i}H\left(\frac{\hat\alpha ^{\top}X_j-\tilde\alpha
^{\top}X_i}{h_2}\right)},
\end{split}
\end{eqnarray}
where $L(\cdot)$ and $H(\cdot)$ are two kernel functions, $h_0$ and $h_2$ are  bandwidths and $\hat \alpha$ is an
estimator derived by a sufficient dimension reduction method
that was described before.

%

\subsection{ Estimation of $\hat \Delta_\tau(x_{10})$}

After having the estimation of $p(\cdot)$, we now proceed to
the step of estimating $\Delta_\tau(x_{10})$. As the estimation procedures for $q_{1,\tau}(x_{10})$ and $
q_{0,\tau}(x_{10}) $ are similar, we only present the detail for $q_{1,\tau}(x_{10})$ and give
the estimator of $q_{0,\tau}(x_{10})$ directly without any
more explanation.

First, we consider the oracle case with the given $p(X)$ and denote the corresponding oracle $CQTE$($OCQTE$) estimator as $\hat \Delta^{ocqte}_\tau(x_{10})$. Note
that for any value $X_{1i}$ that is close to $x_{10}$, Taylor
expansion yields that
$q_{1,\tau}(X_{1i})\approx q_{1,\tau}(x_{10})+q_{1,\tau}'(x_{10})(X_{1i}-x_{10}).$
Thus,  we can use the  minimizer of the following loss
function to define an estimator of ${q}_{1,\tau}(x_{10})$:
\begin{eqnarray}
\begin{split}
 \qquad(\hat{q}_{1,\tau}^{ocqte}(x_{10}),\hat{q'}^{ocqte}_{1,\tau}(x_{10})) =\mathop{\arg\min}_{a,b}
\sum_{i=1}^{n}\frac{D_i}{{p}(X_i)}
\rho_{\tau}(Y_i-a-b(X_{1i}-x_{10}))
K\left(\frac{X_{1i}-x_{10}}{h}\right)\label{equ_sum_1},
\end{split}
\end{eqnarray}
where $K(\cdot)$ is the kernel function and $h$ is the bandwidth.
Similarly, we can define an estimator of
${q}_{0,\tau}(x_{10})$ under the same paradigm:
\begin{eqnarray}\label{equ_sum_0}
\begin{split}
 \qquad (\hat{q}^{ocqte}_{0,\tau}(x_{10}),\hat{q'}^{ocqte}_{0,\tau}(x_{10}))
 =\mathop{\arg\min}_{a,b}
\sum_{i=1}^{n}\frac{1-D_i}{1-{p}(X_i)}
\rho_{\tau}(Y_i-a-b(X_{1i}-x_{10}))K\left(\frac{X_{1i}-x_{10}}{h}\right).
\end{split}
\end{eqnarray}

{ Note that we use the local constant (Nadaraya-Watson (N-W) method) and local linear smoother to estimate $p(\cdot)$ and the function ${q}_{0,\tau}(x_{10})$ respectively. This  is mainly because of the following considerations. We note that the asymptotic bias of ${q}_{0,\tau}(x_{10})$ has no relationship with $\hat p(X)$ as long as its convergence rate can be fast sufficiently. We then  use a simpler estimation for $\hat p(x)$ for ease of exposition,  and the local linear smoother for ${q}_{0,\tau}(x_{10})$ such that the asymptotic analysis can be carefully worked out.
}

Give two
assumptions below. Recall  the definition of high order kernel in the literature. We say a function $g: R^r\rightarrow R$ is a kernel of order $s$ if it integrates to one over $R^r$, and $\int u^{p_1}\cdots u^{p_r}g(u)du=0$ for all nonnegative integers $p_1,\cdots,p_r$ such that $1\leq \sum_{i}p_i<s, $ and it is nonzero when $ \sum_{i}p_i=s. $
{
\begin{itemize}
\item Assumption 3 (on distribution):

(i) The support of the $k$-dimensional covariate vector $X$, $\chi$, is a Cartesian product of compact intervals. The density functions of $X$ and $(X_1,\alpha^{\top}X)$ are bounded away from zero and infinity and $s_1\geq 2$ times continuously differentiable.\\
(ii) The density function of $X_1$, $f(X_1)$,  and the conditional density $f_j(Y(j)\mid X_1)$, are bounded away from zero and infinity and continuously differentiable.
\item Assumption 4 (on kernel function):

 $K(u)$ is a kernel of order $s_1$, is
    symmetric around zero, and is $s^*$ times continuously differentiable.
 \item Assumption 5: $h\rightarrow 0$, $nh^l\rightarrow
    \infty, nh^{2s_1+l+2}\rightarrow 0$.
\end{itemize}

Assumption~3 is commonly used for nonparametric estimation in the literature.   Assumption~4 is for high order kernel. When $l=1$ and $s_1=2,$ Gaussian kernel satisfies this assumption. Further, the value of $s^*$ depends on the estimation procedure to ensure the function smoothness  which will be used in studying the asymptotic behaviors of the estimators.  To be more specific, $s^*\geq 2$ in the case of $PCQTE$, while $s^*\geq s$ and $s^*\geq s_2$ in the case of $NCQTE$ and $SCQTE$ respectively.  Assumption~5 is a condition on the bandwidth selection.} {Obviously, if we assume $nh^{2s_1+l}\rightarrow 0$,  the CQTE estimators can be asymptotically unbiased. However, to better analyze the bandwidth selection rule, we only assume $nh^{2s_1+l+2}\rightarrow 0$.} {

As the benchmark in the latter comparisons, we first present some results about $OCQTE$. Before stating the asymptotic results, define some
important quantities:
\begin{itemize}
\item[(1)] $f_j(y(j)\mid x_{1})$ to be the value of
the conditional density function of $(Y(j)\mid X_1)$, $j=0,1$ at the  point $(Y(j)=y(j), X_1=x_{1})$;
\item[(2)] $\tilde q_{j,\tau}(X_{1i})=q_{j,\tau}(x_{10})+q'_{j,\tau}(x_{10})
    (X_{1i}-x_{10})$, $m_{j,\tau}(X)
=\frac{E[\mathbb I(Y(j)\leq \tilde q_{j,\tau}(X_{1i}))-\tau\mid
X]}{f_j(q_{j,\tau}(x_{10})\mid x_{10})},j=0,1;$

\item[(3)] $\psi(p(X_i),Z_i)=
\frac{D_i}{p(X_i)} \eta_{1,
\tau}(Y_i) -\frac{1-D_i}{1-p(X_i)} \eta_{0,\tau}(Y_i)$, $\sigma_{ocqte}^2(x_{10})=E(\psi^2(p(X),Z)\mid X_1=x_{10}) =E\bigg(\frac{E\big((\mathbb I(Y(1)\leq  q_{1,\tau}(x_{10}))-\tau)^2\mid X\big)}{p(X)f^2_1(q_{1,\tau}(x_{10})\mid x_{10})}+\frac{E\big((\mathbb I(Y(0)\leq  q_{0,\tau}(x_{10}))-\tau)^2\mid X\big)}{(1-p(X))f^2_0(q_{0,\tau}(x_{10})\mid x_{10})} \mid X_1=x_{10}\bigg)$ with $Z_i=(X_i,D_i,Y_i)$;
\item[(4)]$\mu_{s_1}(K)=\int u_1^{p_1}\cdots u_l^{p_l}K(u)du$ for integers $p_1,\cdots, p_l$ such that $\sum_{i=1}^lp_i=s_1$. $\parallel
K\parallel_2^2=\int K^2(u) du$.
\end{itemize}

 \begin{theorem}  For $OCQTE$, when Assumptions 1 through 5 are satisfied,
\begin{eqnarray*}
& &\sqrt{nh^l}\left(\hat \Delta^{ocqte}_{\tau }(x_{10})-
{\Delta_{\tau}
(x_{10})}\right)=-\frac{1}{\sqrt{nh^l}}\frac{1}{f(x_{10})}
\sum_{i=1}^{n}\psi(p(X_i),Z_i)K_i+o_p(1),
\end{eqnarray*}
and
the asymptotic normality   is \begin{eqnarray*}
\sqrt{nh^l}\left(\hat {\Delta}^{ocqte}_\tau (x_{10})-\Delta_\tau
(x_{10})-b_0(x_{10})\right)
\stackrel{\mathcal{D}}{\longrightarrow}N\left(0,\frac{\parallel
K \parallel_2^2\sigma_{ocqte}^2(x_{10})}
{f(x_{10})}\right), \forall x_{10} \in \Omega
\end{eqnarray*}
where the asymptotic bias is $b_0(x_{10})=O_p(\mu_{s_1}(K)h^{s_1})$. Especially, when $s_1=2$, $b_0(x_{10})=\frac{1}{2}\Delta^{''}_\tau
(x_{10})\mu_{2}(K)h^2.$

\end{theorem}
}

 When the propensity
score $p(X)$ is unknown, 
we  use $\hat p(X)$ to replace $p(X)$
 in $(\ref{equ_sum_1})$ and $(\ref{equ_sum_0})$:
\begin{equation}\label{min_hatp}
\begin{split}
&(\hat{q}_{1,\tau}(x_{10}),\hat{q'}_{1,\tau}(x_{10}))=
\mathop{\arg\min}_{a,b}
\sum_{i=1}^{n}\frac{D_i}{{\hat p}(X_i)}
\rho_{\tau}(Y_i-a-b(X_{1i}-x_{10}))
K_i,\\
&(\hat{q}_{0,\tau}(x_{10}),\hat{q'}_{0,\tau}(x_{10}))=\mathop{\arg\min}_{a,b}
\sum_{i=1}^{n}\frac{1-D_i}{{1-\hat p}(X_i)}
\rho_{\tau}(Y_i-a-b(X_{1i}-x_{10}))
K_i.
\end{split}
\end{equation}

For convenience,  denote the estimator $\hat
\Delta_{\tau }(x_{10})$ incorporated with the parametric estimator
$\hat p(X)$ as $PCQTE(x_{10})$, and with the other two nonparametric and semiparametric  estimators
$\hat
p(X)$ separately as $NCQTE(x_{10})$ and $SCQTE(x_{10})$.
As the asymptotic results  vary with the different estimators $
\hat p(X)$ we  present them in the separate subsections.
{For the sake of comparison, all the CQTE estimators are based on the same bandwidth, $h_1$ and kernel function, $K(\cdot)$.}

\subsection{Asymptotic properties of $\hat
\Delta_\tau(x_{10})$ when $p(X)$ is  estimated}

\subsubsection{ PCQTE}

Give the following
assumptions.
{
\begin{itemize}
\item Assumption 6 (Parametric propensity score estimator): The estimator $\hat\beta$ of the propensity score model $\pi(X,\beta)$, $\beta\in\Theta\subset R^d$, $d<\infty$, satisfies $\mathop{\sup}_{X\in \chi}\mid
    \pi(X,\hat \beta)-\pi(X,\beta_0)
    \mid=O_p(n^{-1/2}),$ where $\beta_0\in\Theta$ such that $p(X)=\pi(X,\beta_0)$ for all $X\in \chi.$
\end{itemize}

Assumption~6 is a typical result if we estimate $p(X)=\pi(X,\beta_0)$ by a parametric model like
a logit model or a probit model based on a linear index via the
maximum likelihood method. {The results are stated in the
following theorem, which implies $PCQTE$ is asymptotically equivalent to $OCQTE$.}

\begin{theorem}\label{pcqte_them}
Suppose that Assumptions 1 through 6 are
satisfied for $s_1\geq 2$. Then,   $PCQTE(x_{10})$ has the
asymptotically linear representation as
\begin{eqnarray*}
& &\sqrt{nh^l}\left(\hat \Delta_{\tau }^{pcqte}(x_{10})-
{\Delta_{\tau} (x_{10})}\right)=-\frac{1}{\sqrt{nh^l}}
\frac{1}{f(x_{10})}\sum_{i=1}^{n}\phi_1(p(X_i),Z_i)K_i+o_p(1),
\end{eqnarray*}
and
the asymptotic distribution  is \begin{eqnarray*}
\sqrt{nh^l}\left(\hat {\Delta}^{pcqte}_\tau (x_{10})-\Delta_\tau
(x_{10})-b_1(x_{10})\right)
\stackrel{\mathcal{D}}{\longrightarrow}N\left(0,\frac{\parallel
K \parallel_2^2\sigma_{pcqte}^2(x_{10})}
{f(x_{10})}\right), \forall x_{10} \in \Omega
\end{eqnarray*}
where $\phi_1(p(X),Z)=\psi(p(X),Z)$,  $b_1(x_{10})=O_p(\mu_{s_1}(K)h^{s_1})$. Further,  the asymptotic variance is $
\sigma_{pcqte}^2(x_{10})=\sigma_{ocqte}^2(x_{10}).
$
\end{theorem}
}



\subsubsection{ NCQTE}\label{nonparametrix+setting}

We  make some additional assumptions about kernel functions $
L(\cdot)$ and bandwidths $h$ and $h_0$ to backup the theoretical development.
{
\begin{itemize}
\item Assumption 7:  $L(u)$ is a kernel of order $s\geq k+l$,
    is symmetric around zero, has finite support $[-M,M]^k$, and its $(s+1)$th derivative is continuous. Further,  the density function of $X$, $f_x(X),$ is $s$ times continuously
    differentiable and bounded away from zero and infinity.
\item Assumption 8: $h_0\rightarrow 0$ and
    $\log(n)/(nh_0^{k+s})\rightarrow 0.$

\item Assumption 9: $h_0^{2s}h^{-2s-l}\rightarrow 0$, $nh^{-l}h_0^{2s} \rightarrow 0$.
\end{itemize}

Assumption~7 is also to ensure the smoothness of the density function. Assumption~8 and 9 are the technical conditions to guarantee the existence of the limiting distribution when we need to prove the asymptotic negligibility of all remainder terms. These are because of the involvement of two bandwidths. 

\begin{theorem}\label{ncqte_them}
 Suppose that Assumptions $1$ through $5$ and $7$
through $9$ are satisfied for some $s^*\geq s\geq k+l$, for each point
$x_{10}\in \Omega$, the asymptotically linear representation of
    $NCQTE(x_{10})$ is
 \begin{eqnarray*}
& &\sqrt{nh^l}\left(\hat \Delta^{ncqte}_{\tau }(x_{10})-
{\Delta_{\tau}  (x_{10})}\right)=-\frac{1}{\sqrt{nh^l}}
\frac{1}{f(x_{10})}\sum_{i=1}^{n}\phi_2(p(X_i),Z_i)K_i+o_p(1).
\end{eqnarray*}
 The asymptotic distribution of $\hat \Delta^{ncqte}_\tau(x_{10})$ is
\begin{eqnarray*}
\sqrt{nh^l}\left(\hat {\Delta}^{ncqte}_\tau (x_{10})-\Delta_\tau
(x_{10})-b_2(x_{10})\right)
\stackrel{\mathcal{D}}{\longrightarrow}N\left(0,\frac{\parallel
K \parallel_2^2\sigma_{ncqte}^{*2}(x_{10})}
{f(x_{10})}\right),
\end{eqnarray*}

where $\phi_2(p(X_i),Z_i)=\psi
(p(X_i),Z_i)-n_p(X_i)\epsilon_i$,
$n_p(X_i)=\frac{m_{1,\tau}(X_i)}
{p(X_i)}+\frac{m_{0,\tau}
(X_i)}{(1-p(X_i))},$ $b_2(x_{10})=O_p(\mu_{s_1}(K)h^{s_1})$, $\epsilon=D-p(X)$.
\end{theorem}
}
Rewrite $\phi_2(p(X_i),Z_i)$ as
$
\phi_2(p(X_i),Z_i) = \frac{D_i(\eta_{1,\tau}
(Y_i)-m_{1,\tau}(X_i))}
{p(X_i)}-\frac{(1-D_i)(\eta_{0,\tau}(Y_i)-m_{0,\tau}(X_i))}{1-p(X_i)}
 +m_{1,\tau}(X_i)-m_{0,\tau}(X_i).$
We then get the following corollary.
\begin{corollary}\label{cor_nop} Under the regularity conditions in the previous theorems,
   \begin{eqnarray*}
   (1)  \sigma_{ncqte}^{2*}(x_{10})&=&E\left[
   (m_{1,\tau}(X)-m_{0,\tau}(X))^2+\frac{\sigma_{\tau,1}^2(X)}
{p(X)}+\frac{\sigma_{\tau,0}^2(X)}{1-p(X)}\mid X_1=x_{10}
 \right],\\
 (2)
 \sigma_{pcqte}^2(x_{10})&=&\sigma_{ncqte}^{*2}(x_{10})+E\left[p(X)(1-p(X))
 \left(\frac{m_{1,\tau}(X)}{p(X)}+\frac{m_{0,\tau}(X)}{1-p(X)}
 \right)^2  \mid
 X_1=x_{10}\right]\geq \sigma_{ncqte}^{*2}(x_{10}),
 \end{eqnarray*}
 where $\sigma_{\tau,j}(x)=Var\bigg(\frac{\mathbb I(Y(j)\leq  q_{j,\tau}(X_{10}))-\tau}{f_j(q_{j,\tau}(x_{10})\mid x_{10})}\mid X\bigg).$
\end{corollary}

  \begin{rema}  {{\bf $Corollary$ $2.1$} implies that for any $x_{10}$, $NCQTE \preceq PCQTE \cong OCQTE$.
}{As discussed before, we can not show whether $NCQTE(X_1)$ is the most efficient $CATE$ estimator as the standard semiparametric efficient theory is invalid for an unknown function. This phenomenon is very different from the unconditional quantities, e.g. $ATE$ and $QTE$, which can  achieve the semiparametric efficient bound when nonparametrically estimated propensity score is used.
But looking at  all the asymptotic variance functions, $NCQTE(x_{10})$  is the most efficient estimator and thus, we conjecture that $NCQTE(x_{10})$ would achieve an efficient bound in certain sense. This deserves a further study.}
\end{rema}
\subsubsection{ SCQTE}
{
If we  postulate that  the information about $D$ from $X$ can be
completely captured by $q$ linear combinations $\alpha^{\top}X$ of $X$ with $l\leq q\ll k$, we can then estimate the propensity
score function $p(X)=p(\alpha^\top X)$  with
$\alpha^{\top}X$ rather than the original $X$ to avoid the curse of dimensionality. To this end,  we
can use a lower dimensional kernel function $H(u)$, instead of a
high dimensional kernel function $L(X)$ to get the local smooth
estimator $\hat p(\alpha^\top X)$ of $p(X)$.}

%


 we first correspondingly rectify the assumptions related
 to $L(X)$ and bandwidth in Subsection~ \ref{nonparametrix+setting}.

\begin{itemize}
\item Assumption 7': $H(u)$ is symmetric around zero,  has finite support $[-M,M]^{q}$, and
    is $s_2\geq q+l$ times continuously differentiable. The density function of $\alpha^{\top}X$, $f_\alpha(\alpha^{\top}X)$ is $s_2$ times continuously
    differentiable.
\item Assumption 8': $h_2\rightarrow 0$ and
    $\log(n)/(nh_2^{s_2+q})\rightarrow 0.$

\item Assumption 9': $h_2^{2s_2}h^{-2s_2-l}\rightarrow 0$, $nh^{-l}h_2^{2s_2} \rightarrow 0$.

\item Assumption 10': $\hat \alpha$ is a $root$-$n$ consistent estimator of $\alpha$ and $q$, the dimension of $\alpha$, is given with $l\leq q\ll k$.
\end{itemize}
{
Since the propensity score of $SCQTE(x_{10})$ is based on $\alpha^{\top}X$,
Assumptions 7' through 9' are adjusted to those in the case of $NCQTE(x_{10})$ and play the same role.
}
{
Further, we define some nations for ease of interpretation. Let $card(A)$ be  the cardinality of set $A$. $A\subseteq B$ means $A\cap B=A$, that is, all elements of $A$
    are also elements of $B$, while $A\subsetneq B$ means $A\subseteq B$ but $card(A)<card(B)$.
     $A\sqsubset^t B$ means $A\cap B=C$ with the cardinality $card(C)=t$, that is, there
    only exist $t$ elements of $A$ belonging to $B$. Especially when $t=0$,
    it means none of elements in the set $A$ are the elements in set $B$.
}

\begin{theorem}{\label{sem_them}}
Suppose the assumptions 1 through 4 and 7'
through 9' are satisfied for $s^*\geq s_2\geq q+l$,
  the following statements  hold for each point
$x_{10}\in \Omega$:
\begin{itemize}
\item [(1)] When $X_1 \sqsubset^{l-r} \alpha^{\top}X$ with $0< r\leq l$ and $s_2(2-l/r)+l>0$,
 the asymptotically linear representation of
    $SCQTE(x_{10})$ is
 \begin{eqnarray*}
& &\sqrt{nh^l}\left(\hat \Delta^{scqte}_{\tau }(x_{10})-
{\Delta_{\tau} (x_{10})}\right)=-\frac{1}{\sqrt{nh^l}}
\frac{1}{f(x_{10})}\sum_{i=1}^{n}\phi_3(p(\alpha^\top
X_i),Z_i)K_i+o_p(1)
\end{eqnarray*}
and the asymptotic distribution of $\hat \Delta^{scqte}_\tau(x_{10})$ is
\begin{eqnarray*}
\sqrt{nh^l}\left(\hat {\Delta}^{scqte}_\tau (x_{10})-\Delta_\tau
(x_{10})-b_3(x_{10})\right)
\stackrel{\mathcal{D}}{\longrightarrow}N\left(0,\frac{\parallel
K \parallel_2^2\sigma_{scqte}^2(x_{10})}
{f(x_{10})}\right),
\end{eqnarray*}
\item[(2)]
 When $X_1 \subsetneq \alpha^{\top}X$,
 the asymptotically linear representation of
    $SCQTE(x_{10})$ is
 \begin{eqnarray*}
& &\sqrt{nh^l}\left(\hat \Delta^{scqte}_{\tau }(x_{10})-
{\Delta_{\tau} (x_{10})}\right)=-\frac{1}{\sqrt{nh^l}}
\frac{1}{f(x_{10})}\sum_{i=1}^{n}\phi_3^*(p(\alpha^\top
X_i),Z_i)K_i+o_p(1),
\end{eqnarray*}
and the asymptotic normality of $\hat \Delta^{scqte}_\tau(x_{10})$ is
\begin{eqnarray*}
\sqrt{nh^l}\left(\hat {\Delta}^{scqte}_\tau (x_{10})-\Delta_\tau
(x_{10})-b_3^*(x_{10})\right)
\stackrel{\mathcal{D}}{\longrightarrow}N\left(0,\frac{\parallel
K \parallel_2^2\sigma_{scqte}^{*2}(x_{10})}
{f(x_{10})}\right),
\end{eqnarray*}
\item [(3)] $SCQTE(x_{10})$ has a limiting variance that is smaller or equal to those of  $PCQTE(x_{10})$ and the oracle  $CQTE(x_{10})$ as below: 
\begin{eqnarray*}
&&\sigma_{pcqte}^2(x_{10})=\sigma_{ocqte}^2(x_{10})=\sigma_{scqte}^2(x_{10})=\sigma_{scqte}^{*2}(x_{10})\\
&&+ E\left[p(\alpha^{\top}X)(1-p(\alpha^{\top}X))\left(\frac{m_{1,\tau}(\alpha^{\top}X)}{p(\alpha^{\top}X)}
 +\frac{m_{0,\tau}(\alpha^{\top}X)}{1-p(\alpha^{\top}X)}\right)^2\mid
X_1=
x_{10}\right]\geq \sigma_{scqte}^{*2}(x_{10}).
\end{eqnarray*}
\end{itemize}

where $\phi_3(p(\alpha^\top X_i),Z_i)=\psi
(p(X_i),Z_i)$, $b_3(x_{10})=O_p(\mu_{s_1}(K)h^{s_1})$.
$\phi^*_3(p(\alpha^{\top}X_i),Z_i)=\psi(p
(\alpha^{\top}X_i,Z_i))-e_p(\alpha^{\top}X_i)\epsilon_i,$
$e_p(\alpha^\top x_i)=\frac{m_{1,\tau}(\alpha^\top X_i)}
{p(\alpha^\top X)}+\frac{m_{0,\tau}
(\alpha^\top X_i)}{1-p(\alpha^\top X)},$ $b_3^*(x_{10})=O_p(\mu_{s_1}(K)h^{s_1})$,
and $\epsilon=D-p(\alpha^\top X)$.
\end{theorem}
{
\begin{rema}
We should also note that when $\alpha^\top X=X_1$, we have $E(m_{j,\tau}(\alpha^\top X)\mid X_1=x_{10})=0, j=0,1$, and the asymptotic variance as $$\sigma_{pcqte}^2(x_{10})=\sigma_{scqte}^{*2}(x_{10})
=\frac{\tau(1-\tau)}{p(x_{10})f^2_1(q_{1,\tau}(x_{10})\mid x_{10})}+\frac{\tau(1-
\tau)}{(1-p(x_{10}))f^2_0(q_{0,\tau}(x_{10})\mid x_{10})}.$$
That implies, when $\alpha^\top X=X_1$, $SCQTE$ cannot be more efficient than $PCQTE$ even when the propensity score is estimated nonparametrically. This is an essential difference from the unconditional counterpart. But when $X_1 \subsetneq\alpha^{\top}X$, the nonparametric structure of $SCQTE(x_{10})$ estimator does play a positive role in efficiency. In Subsection~\ref{add_section} below, we give some more discussions and more general results to provide a relatively complete picture of estimation efficiency in this field. 
\end{rema}

{
\subsection {Further studies about the role of  propensity score in efficiency   }\label{add_section}

The above results about $CQTE$ estimators present two interesting phenomena. In the scenario Theorem~\ref{ncqte_them} presents, $NCQTE(x_{10})$ can be asymptotically more efficient than $PCQTE(x_{10})$ and the oracle $CQTE(x_{10})$. Yet, in the scenario Theorem~\ref{sem_them} designs, $SCQTE(x_{10})$ cannot always be so  although $SCQTE(x_{10})$ also uses the nonparametric method to estimate the  propensity score.
{
This motivates us to further investigate the role of the estimated propensity score in the asymptotic behaviors of the $CQTE$ estimators. At first glance, it seems that the different asymptotic behaviors  are because of different estimation methods for propensity score. Comparing $SCQTE(x_{10})$ with  $NCQTE(x_{10})$ and  the technical proofs for Theorem~~\ref{ncqte_them} with that for Theorem~\ref{sem_them} in Appendix,  we note that there are two factors playing the important role in the estimation efficiency: how fast is the convergence rate of the estimated propensity score and whether  $X_1$ is  a strict subset of the true arguments of the propensity score. We then separately discuss them.
 }
\subsubsection{The role of convergence rate of the estimated propensity score }
{
From the technical proofs and the main differences between $PCQTE$  and $NCQTE$ / $SCQTE$  we can see that fast rate of convergence can make the first order expansion of the estimated propensity score such as $PCQTE$ vanish while slow rate such as for $NCQTE$ / $SCQTE$ cannot cancel off it and thus enhance the asymptotic efficiency due to a negative correlation with the leading term. Thus, 
 $NCQTE$ and $SCQTE$  can lose their efficiency superiority if the nonparametric estimations converge faster with higher order smoothness as concluded in the following corollary.
}

 \begin{corollary}\label{add1}
    In addition to the conditions in Theorem~\ref{ncqte_them} or Theorem~\ref{sem_them} respectively with replacing the assumptions on the bandwidths $h_p$ and $h$ by $\sqrt {nh^l}(h_p^s+\sqrt {\log (n)/nh_p^k})=o(1)$ for some number $s$, $NCQTE(x_{10})$ and $SCQTE(x_{10})$ have the same asymptotic distribution as $PCQTE(x_{10})$. That is,
     $$OCQTE(x_{10})\cong PCQTE(x_{10})\cong NCQTE(x_{10}) \cong SCQTE(x_{10}).$$
 \end{corollary}

\begin{rema} Obviously,  the above discussion is
mainly for theoretical investigation. In practice, it makes no sense to choose such bandwidths to use. But the discussion is still helpful for us to better understand the estimation mechanisms.  These results further reveal the essential differences between the conditional and unconditional structure. As we known under  the unconditional structure, the estimator requires a standardizing constant $\sqrt n$ to derive its limiting distribution. Any  estimator of the propensity score is at the  rate of order $1/\sqrt n$ or slower, the impact from the estimated propensity score will play role in reducing the asymptotic variance. {See \cite{hir2003} and \cite{liu2018}.} In contrast, under the conditional structure,  the estimator requires a standardizing constant $\sqrt {nh^l}$ to approximate its limiting distribution. When an  estimator of the propensity score is at a rate $b_n$ faster than $1/\sqrt{nh^l}$ such as the case where the propensity score is parametric with the rate of order $1/\sqrt n$, the impact from the estimated propensity score will play no role in the asymptotic variance reduction.
{For CATE, \cite{abrevaya2015} showed the case with the parametric propensity score, but did not include the discussion on nonparametric and semiparametric cases.} Thus, for the estimated semiparametric and nonparametric propensity score, delicately choosing the bandwidth to obtain a proper rate of convergence becomes vital for the estimation efficiency.
It is clear that in the above corollary, the condition $\sqrt {nh^l}(h_p^s+\sqrt {\log (n)/nh_p^k})=o(1)$ is much stronger than the assumptions  in Theorems~2.3 and 2.4, but is still possible to choose such bandwidths as long as the involved functions are sufficiently smooth and high order kernels are used. This is because for the nonparametric estimation, the rate of convergence can be as close to $1/\sqrt n$ as possible when the function is very smooth. {However, utilizing a high order kernel for regression fit means we would assign negative weights to some range of the data, which can be an undesirable side-effect. See \cite{li2007}. }
\end{rema}
{ \begin{rema}The theorems in this paper also add new insights about the super-efficiency phenomenon found in missing data and treatment effect area. That is, for unconditional treatment effect, generally inverse of propensity score-based estimators with estimated propensity score is more efficient than the one with true propensity score. As discussed above, estimating propensity score is not necessary to play role in the asymptotic variance reduction. 
\end{rema}}

\subsubsection{The effect of the affiliation of $X_1$ to the  set of true arguments of propensity score function}
{
As pointed out before, $CQTE$ is a function of $X_1$, and its affiliation to the set of all arguments of the propensity score plays role for estimation efficiency. Recall that in the scenario Theorem~\ref{sem_them} discusses, the asymptotic distribution of $SCQTE(x_{10})$ depends on the relationship between $X_1$ and $\alpha^{\top}X$. Therefore, under the  constraint~2 below, it can be expected that $NCQTE(x_{10})$ will have similar properties, namely, the affiliation of $X_1$ to $\tilde X$ (or $\alpha^{\top}X$) should affect the asymptotic distribution of $NCQTE(x_{10})$ (or $SCQTE(x_{10})$).
\begin{eqnarray}\label{cons_2}
\mbox{Constraint~2}:\ \ D\perp X \mid \tilde X.
\end{eqnarray}
Thus { we call $\tilde X$ is the set of true arguments of propensity score and $p(X)=p(\tilde X).$} Obviously, when $\tilde X=X$, $X_1\subseteq \tilde X$, $NCQTE(X_{10})$ can be more efficient than $PCQTE(X_{10})$ by Theorem~\ref{ncqte_them}. Let $X_1 \sqsubset^{l-r} \tilde X$ mean $X_1$ is $l-r$ components of $\tilde X$. We will see that when $X_1 \sqsubset^{l-r} \tilde X\subsetneq X$, the situation will be different as concluded by the following corollary.
\begin{corollary}\label{add2}
  Suppose that there is a given $\tilde X$ such that  $D\perp X\mid \tilde X$ with $X_1 \sqsubset^{l-r} \tilde X\subsetneq X$ and $0< r\leq l$.  Then if the propensity score
$p(\tilde X)$ is estimated by basing on $\tilde X$ rather than  $X$, under the conditions in Theorem~\ref{ncqte_them} and $s(2-l/r)+l>0$, $NCQTE(x_{10})$ has the same asymptotic distribution as $PCQTE(x_{10})$. Then
 $ NCQTE \cong PCQTE.$
 \end{corollary}
}


%

{
Further, we clarify the relation between  $NCQTE(X_{10})$ and $SCQTE(X_{10})$ when both $X_1 \subsetneq \tilde X=X$ and $X_1 \subsetneq
     \alpha^{\top}X$ hold.

 \begin{corollary}\label{var_com}
 Suppose all the assumptions listed above and the two assumptions~(\ref{cons_1}) and (\ref{cons_2}) are satisfied, namely $p(X)=p(\tilde X)=p(\alpha^{\top}X)$, and  $X_1 \subsetneq \tilde X=X$ and $X_1 \subsetneq
     \alpha^{\top}X$, we have the following asymptotic variance functions of $SCQTE$ and $NCQTE$:\\
$$\sigma_{scqte}^{*2}(z)=\sigma_{ncqte}^{*2}(z)+
E\bigg[p(\alpha^{\top}X)(1-p(\alpha^{\top}X))\bigg\{\frac{\Delta m_{1,\tau}}{p(\alpha^{\top}X)}
 +\frac{\Delta m_{0,\tau}}{1-p(\alpha^{\top}X)}\bigg\}^2\mid X_1=x_{10}\bigg],$$
where $\Delta m_{j,\tau}=m_{j,\tau}(X)-m_{j,\tau}(\alpha^{\top}X).$
%
%
%
  \end{corollary}

We are now in the position to summarize all results about the affiliation effect of $X_1$.
\begin{itemize}
\item[(1)] $NCQTE(X_{10}) \preceq SCQTE(X_{10})\preceq PCQTE(X_{10})\cong OCQTE(X_{10})$, for $X_1 \subsetneq \tilde X=X$ and $X_1 \subsetneq
     \alpha^{\top}X$;
\item[(2)]  $NCQTE(X_{10}) \preceq SCQTE(X_{10})\cong PCQTE(X_{10})\cong OCQTE(X_{10})$, for $X_1 \subsetneq \tilde X=X$ and $X_1 \sqsubset^{l-r}
     \alpha^{\top}X$;
\item[(3)] $NCQTE(X_{10}) \cong SCQTE(X_{10})\cong PCQTE(X_{10})\cong OCQTE(X_{10})$, for $X_1 \sqsubset^{l-r}\tilde X\subsetneq X$ and $X_1 \sqsubset^{l-r}
     \alpha^{\top}X$.
\end{itemize}

\begin{rema} It is very interesting that whether $NCQTE(X_{10})$ and $SCQTE(X_{10})$ can be asymptotically more efficient also relies on whether the given covariates are a strict subset of the arguments of {the  propensity score.}  This important observation is not easy to explain, but might be because of the following. Note that  when it does not include the given covariates, then under the conditional structure, the estimated propensity score is conditionally independent of the conditional treatment effect and then plays little role for the asymptotic property of the estimated treatment effect. It deserves a further study to confirm this explanation.
\end{rema}%
}
\vspace{-10pt}

{\subsection{ An extension to the large $k$ setting}

When the dimension $k$ of $X$ is large or even larger than the sample size $n$, the   CQTE estimation needs a further dimension reduction combining variable selection and  a post-selection estimation. Then the relevant investigation can be conducted. To this end, we can modify the estimation procedure of $PCQTE$ and $SCQTE$ when there is a sparsity structure of propensity score model. That is, only a relatively small number of important covariates are selected to treatment assignments while the rest are treated as unimportant ones. We give some descriptions on the basic ideas below.

More specifically, for $PCQTE$, we can replace the propensity score $\hat p(X_i)=\pi(X_i,\hat \beta)$ by a penalized maximum likelihood estimator $\pi(X_i,\hat \beta_c)$ (e.g. \citet{fan2004} ) where $\hat \beta_c$ is obtained by maximizing
\begin{eqnarray}
 \sum_{i=1}^n(D_i \log
\pi(X_i,\beta)+(1-D_i)(1-\log \pi(X_i,\beta) )-\sum_{j=1}^kR_{\lambda}(\beta_j).
\end{eqnarray}
Here $R_{\lambda}(\beta_j)$ is a penalized function designed to select important variables with the regularization parameter $\lambda$ being chosen by cross-validation.
There are several choices of $R_{\lambda}(\beta_j)$, such as the Lasso \citep{tibshirani1996} and smoothly clipped absolute deviation(SCAD) penalty \citep{fan2001}. We can then obtain $PCQTE$ based on $\pi(X_i,\hat \beta_c)$ by solving the optimal problem (\ref{min_hatp}).

For $SCQTE$, we can also replace the classical dimension reduction method with a sparse dimension reduction method \citep{wang2018}, which combines variable selection and model-free sufficient dimension reduction together, to estimate the propensity score. A much relevant literature is \citet{ma2019}, who also proposed a new sparse dimension reduction method to estimate propensity score  for estimating the average treatment effect.

}{
\subsection{Estimation for asymptotic variance }

{We also very briefly describe the issue of estimating the asymptotic variance functions.  In the following, we take $PCQTE$ as an example to briefly describe an estimation procedure,  the variance functions of the other $CQTE$ estimators can be similarly estimated.

Recall that the asymptotic variance of $PCQTE$ is $\frac{\parallel
K \parallel_2^2\sigma_{pcqte}^2(x_{10})}
{f(x_{10})}$, we then need to consistently estimate $\sigma_{pcqte}^2(x_{10})$ and $f(x_{10})$. For $f(x_{10})$, the nonparametric kernel estimation, $\frac{1}{nh^l}\sum_{i=1}^nK\left(\frac{X_{1i-x_{10}}}{h}\right)$ can be used. For $\sigma_{pcqte}^2$, the kernel estimator is as
 \begin{eqnarray}
 \hat\sigma_{pcqte}^{2}(x_{10})=\bigg[ \frac{1}{nh^l} \sum_{i=1}^n\{\hat \phi_1(\pi(X_i,\hat \beta),Z_i)\}^2 K\left(\frac{X_{1i}-x_{10}}{h}\right) \bigg]/\hat f(x_{10}).
 \end{eqnarray}
Here $\hat \phi_1(\hat p(X_i),Z_i)=\frac{D_i\hat\eta_{1,\tau}
(Y_i)}
{ \pi(X_i,\hat \beta)}-\frac{(1-D_i)\hat \eta_{0,\tau}(Y_i)}{1- \pi(X_i,\hat \beta)}$ and $\hat \eta_{j,\tau}(Y_i)=\frac{\mathbb I(Y_i\leq \hat q_{j,\tau}(x_{10}))}{\hat f_j(\hat q_{j,\tau}(x_{10})\mid  x_{10} )},\ j=0,1$, and $\hat f_j(\hat q_{j,\tau}(x_{10})\mid  x_{10} )=\frac{\sum_{i=1}^n\hat w_jK\left[(Y_i-\hat q_{j,\tau}(x_{10}))/h\right]K\left[(X_{1i}-x_{10})/h\right]}
{\sum_{i=1}^nK\left[(X_{1i}-x_{10})/h\right]}$ with $\hat w_1=\frac{D_i}{\pi(X_i,\hat \beta)}$ and $\hat w_0=\frac{1-D_i}{1-\pi(X_i,\hat \beta)}$. As all are related to nonparametric kernel estimations, the consistency can also be expected. However, we also see that it involves many unknowns, the estimation may not be efficient sufficiently in finite sample scenarios.

An alternative  is the nonparametric bootstrap approximation \citep{efron1979}, which is often useful in practice. The procedure can be described by the following steps: given $X_1=x_{10}\in \Omega$,
\begin{itemize}
\item Step~1: Given original random sample $\{(Y_i,X_i,D_i):i=1,\cdots,n\}$, obtain the maximum likelihood propensity score estimator $\pi(X_i,\hat \beta)$ and $\hat \Delta_{\tau}^{pcqte}(x_{10})$ as described before;
\item Step~2: Generating the $b$-th bootstrapped sample $\{(Y^{b}_i,X^{b}_i,D^{b}_i):i=1,\cdots,n\}$, $b=1,\cdots,B$ with replacement from  $\{(Y_i,X_i,D_i):i=1,\cdots,n\}$. For each bootstrapped sample, compute $\hat \pi(X_i,\hat \beta_b )$ and $\hat \Delta_{\tau,b}^{pcqte}(x_{10})$;
\item Step~3:  The estimator of the asymptotic variance of $\hat \tau_0(x_{10})$ can be obtained by the empirical variance of $( \hat \Delta_{\tau,1}^{pcqte}(x_{10}),\cdots,  \hat \Delta_{\tau,B}^{pcqte}(x_{10}))$:
\begin{eqnarray}
\hat {Var}[\hat \Delta_{\tau}^{pcqte}(x_{10})]=\frac{1}{B-1}\sum_{b=1}^B
\left[ \hat \Delta_{\tau,b}^{pcqte}(x_{10})-\hat \Delta_{\tau}^{pcqte}(x_{10})\right]^2.
\end{eqnarray}
\end{itemize}

As this is not the focus of this paper, we then do not give more details about their asymptotic properties.
}

}
\subsection{The bandwidth selection rule}

Note that
$PCQTE$ only involves one bandwidth  $h$  used in the integration step. Minimizing the asymptotic MISE under $s_1=2$ leads to the asymptotically optimal bandwidth as $
h^p_{opt}=\left(\frac{\|K\|_2^2\int                                                                (\sigma^2_1(X_1)/f(X_1))dX_1}{\mu^2_2(K)\int (\Delta''_\tau(X_1))^2dX_1}\right)n^{-1/5}:=C_\tau n^{-1/5}.
$
However,  the bandwidth selection is always very critical for  the asymptotic behaviors of $SCQTE$ and $NCQTE$. We need to delicately choose the bandwidths $h_0$ and $h_2$ that are used  in the estimated propensity score separately for $NCQTE$ and $SCQTE$, and $h$ that is for the final estimator. Since the bandwidth selection procedure is very much complicated to balance the magnitudes between these bandwidths,  $SCQTE$ and $NCQTE$ can be sensitive to the selected bandwidths. The results in simulation show this phenomenon. Further, we should note that the first-order asymptotic theory derived in this paper cannot provide the idea on the optimal bandwidth selection.


{
Thus, we turn to use the rule of thumb to guide the bandwidth selection. Take the case of $NCQTE$ as an example. Recall that the corresponding bandwidths, which are respectively denoted as
$h_0, h,$ should satisfy the Assumptions $5,7-9$.
Thus to fulfill the assumptions, take
\begin{equation}\label{band_rule}
 h=a\cdot n^{\frac{-1}{l+2s_1-\delta}},  h_0=a_1\cdot
n^{\frac{-1}{k+s+\delta_0}}, \quad \mbox{ for } \quad a>0,\delta>0,
  a_1>0,\delta_0>0.
\end{equation}
Note that $\delta$ and $\delta_0$ can be made as small as necessary or desired, thus we set them as zero. Further the order of kernel function are set as:
 $s_1=s$ and $s=k+l+\delta^*$. Here when $k+l$ is even, $\delta^*=1$, otherwise $\delta^*=0$. When it comes to $SCQTE$,
let \begin{equation}\label{band_rule2}
 h=b\cdot n^{\frac{-1}{l+2s_1-\delta}},  h_2=b_1\cdot
n^{\frac{-1}{q+s_2+\delta_0}}, \quad \mbox{ for } \quad b>0,\delta>0, b_1>0,\delta_0>0.
 \end{equation}
 Obviously, we just need to replace the role of $k$ in $NCQTE$ by $q$ . The  above rule is not the unique way, but  is easy to implement, and  thus is a good way in practice.}

%
%
\section{Simulation studies}
\subsection{Preliminary of the simulation}\renewcommand{\theequation}{3.\arabic{equation}}
\setcounter{equation}{0}

In this section, we aim to compare the finite sample
performance of the proposed estimators, taking $OCQTE$ as the benchmark to examine the aforementioned theoretical results. For ease of exposition, we only consider the case of $X_1\in R$, i.e. $l=1$.
{
Further, to better analyze the performance of $NCQTE$, we only
consider $k=dim(X)\in \{2,4\}$, which turns out that this setup is sufficiently informative to show $NCQTE$ obtains the efficiency superiority when $k=2$ and loses this superiority when $k=4$, due to the dimensionality problem in nonparametric estimation.
}

Consider the following heteroscedasticity models for $k=2$ and $k=4$ respectively:
\begin{description}
\item [\bf{Model 1:}]
    $Y(0)=0$, and
    $Y(1)=X_1+X_2+|X_1|\epsilon_1$, $p_1(X)=\frac{\exp(\alpha_1^{\top}X)}
{1+\exp(\alpha_1^{\top}X)}, $\label{model_1}
\item [\bf{Model 2:}]$Y(0)=0$, and
    $Y(1)=X_1+X_2+|X_1|\epsilon_1$, $p_2(X)=\frac{\exp(\alpha_2^{\top}X)}
{1+\exp(\alpha_2^{\top}X)}, $\label{model_1}
\item[\bf{Model 3:}]$Y(0)=0$, and
    $Y(1)=X_1+X_2+X_3+X_4+|X_1|\epsilon_1$, $p_3(X)=\frac{\exp(X_1^2+\alpha_3^{\top}X)}
{1+\exp(X_1^2+\alpha_3^{\top}X)} .\label{model_2}$
\end{description}
 Here $\alpha_1=(1,1)^{\top}, \alpha_2=(0,1)^{\top}$,  $\alpha_3=(0,1/\sqrt{3},1/\sqrt{3},1/\sqrt{3})^{\top}$.

Obviously, under $Model~1$, when $p(X)=p_1(X)$, $X_1\subseteq \tilde X=X$ but $X_1 \sqsubset^{0} \alpha_1^{\top}X$, which is designed to examine whether $NCQTE$ is the most efficient estimator while $SCQTE$ is asymptotically similar to $PCQTE$ and $OCQTE$. As for $p_2(X)$ in $Model~2$, since $X_1 \sqsubset^{0} \tilde X=X_2$ and $X_1  \sqsubset^{0} \alpha_2^{\top}X$, it can be expected that all $CQTE$ are asymptotically similar.  $p_3(X)$ in $Model~3$ is set to verify that $SCQTE$ can be more efficient that $PCQTE$. In this propensity score model, $D\perp X\mid \alpha^{\top}X$ with  $\alpha^{\top}=\left(
                 \begin{array}{cccc}
                   1 & 0 & 0 & 0\\
                   0 & 1/\sqrt{3}& 1/\sqrt{3}& 1/\sqrt{3} \\
                 \end{array}
               \right)$.
Here, we use the aforementioned sufficient dimension reduction method, MAVE, to estimate the index $\alpha^{\top}X$.

When $k=2$,  generate $X=(X_1,X_2)^{\top}$,  $\epsilon_i$ from
 $ X_1\sim U(-0.5,0.5),X_2=(1+X_1^2)+\varepsilon_1, $ $\varepsilon_1\sim N(0,0.25^2)$ and $\epsilon_i\sim N(0,1), i=1,2.$
When the dimension of $X$ is $k=4$, the generation of  $X=(X_1,X_2,X_3,X_4)^{\top}$ and
$\epsilon_3$ is
$X_1\sim U(-0.5,0.5), X_2= (1+X_1^2)+\varepsilon_1, X_3= X_1(1+X_1)+\varepsilon_2, X_4=\exp(-1-X_1)+\varepsilon_3,$
where $\epsilon_3 \sim N(0,1)$, $\varepsilon_j\sim N(0,0.25^2)$, $j=1,2,3$.

Thus the corresponding $\tau th$ $CQTE$ $\Delta_\tau(X_1)$ under
the mentioned models are:
\begin{description}
\item [\bf{Model 1\& Model 2:}] $\Delta_\tau(X_1)=X_1+ (1+X_1^2)+F_1^{-1}(\tau),$
\item[\bf{Model~3:}]$\Delta_\tau(X_1)=X_1+(1+X_1^2)+X_1(1+X_1)+\exp(-1-X_1)+F_2^{-1}(\tau),$
\end{description}
where $F_1^{-1}(\tau),F_2^{-1}(\tau)$ are respectively the $\tau$th quantile of $N(0,X_1^2+0.25^2)$ and $N(0,X_1^2+3\times0.25^2).$ Note that $CATE(X_1)$ cannot capture the heteroscedasticity structure of error term while $QTE$ is just a quantity. 
%
%

In the  procedures described here,  we give the estimators of
$\Delta_\tau(x_1) $ at $x_1\in\{-0.2,0,0.2\}$. To save space, we only report the quantile level $\tau=0.5$ in the simulations, while the other quantile levels share similar finite sample performance in the comparisons among the estimators.  Two sample sizes are considered: $n=500$ and $n=1000.$ The replication time is $1500$.  We choose a Gaussian kernel and then high order kernels are derived from it throughout this section.
 Next, the values of the order and smoothness $s_1,s_2,s$ of the kernels and the bandwidths $h,h_0,h_2$ are chosen  via the rules in  (\ref{band_rule}) and  (\ref{band_rule2}). To better examine the performances fairly, the parameters $s_1,h$ for $K(u)$
are the same for all  four $CQTE$ estimators. Since the conditions for $NCQTE$ are more restrictive than those for the other estimators, we first select all parameters for $NCQTE$, i.e. the tuning parameters in  (\ref{band_rule}) and  (\ref{band_rule2}), where we set $a=b$. That means we just need to confirm the values of turning parameters $\{a,a_1,b_1\}$ about the bandwidths. By the rule of thumb, we try two groups of values of $a,a_1,b_1$ to see which ones could make stable performances of estimations, that are, $
Group\ 1: \{a=0.5,a_1=1.1,b_1=1.2\},\  Group\ 2:~
\{a =0.5,a_1=0.9,b_1=1.1\}.$

Further, we should point out that the estimated propensity score is
trimmed to lie in the interval $[0.005,0.995]$.
Bias, standard deviation (SD) and mean squared error (MSE) of
$\hat \Delta_\tau(x_1) $ are used to evaluate the performance of the
involved estimators.

\subsection{Simulation results}

For space saving and better illustration, we in this section only display, when the SDs are used, the asymptotic relative efficiency(ARE) of the estimators against $OCQTE$ under $Group~1$ values of $\{a,\ a_1,\ b_1\}$ in Figure~\ref{group1} to visualize their performances. All simulation results  are presented in Table~1.
Analyzing the simulation results reported in Figure~\ref{group1} and Table~1, we summarize the conclusions as follows.

}

{\it The effect of sample size.} Comparing the estimation effect with different sample sizes in the same model, we can see that a larger sample size leads to  smaller  MSE and $SD$. Across all the models, both the MSE's and $SD$'s of these four estimators with the  sample size $n=500$ are roughly 1.5 times larger than those with $n=1000$.

{\it The effect of dimensionality.}  Even though we consider relatively low dimensions in the simulation settings, the  influence by dimension on the evaluation indexes is still observable. For example, when the dimension $k$ increases from $2$ up to $4$, the results in the table show that $SD$ and $MSE$ obviously increase. But we also point out that, even when $k=4$, the values of  $SD$ and $MSE$ are still small. Further, comparing Model~1 in the $k=2$ setting with Model~3 in the  $k=4$ setting in Figure~\ref{group1}, we can see that when $k=2$, $NCQTE$ is almost uniformly more efficient than the other $CQTE$ estimators even when $n=500$. This is consistent with the asymptotic results in Theorem~$2.3$. However, when $k=4$ as illustrated in Figure~\ref{group1}, $NCQTE$ sometimes loses its efficiency superiority even when $n=1000$. This would be due to the estimation inaccuracy when the dimension is high.

{\it The effect of estimation method.} In the simulations, we compare all four estimators.
In terms of all the evaluation indexes,  $PCQTE$ has, in most cases, similar performance to  $OCQTE$. This well coincides with  the asymptotic properties presented in Theorem~$3.2$. Further, as discussed before, the performances of $NCQTE$ and $SCQTE$ are related to whether the given $X_1$ is included in the set of arguments of the propensity score.  Specifically, when $p(X)=p_1(X)$, i.e. under Model~1,  Figure~\ref{group1} shows that $NCQTE$ is uniformly the most efficient one. While when $p(X)=p_2(X)$, i.e. under Model~2, as $X_1$ is not fully included in the set of  the arguments of $p_2(X)$, Figure~\ref{group1} shows that the estimation efficiency of $NCQTE$ loses.  When $p(x)=p_3(x)$, Figure~\ref{group1} shows that both $NCQTE$ and $SCQTE$ are generally more efficient than $OCQTE$ and $PCQTE$. Further, when $k$ is larger, $SCQTE$ can be, in some cases, more efficient than $NCQTE$. That is, $NCQTE$ is no longer always superior to $SCQTE$. This is mainly because of the dimension reduction structure in $SCQTE$ and thus, less estimation inaccuracy. Thus, the performance of $SCQTE$ could be  more robust than $NCQTE$ against dimensionality.


In summary,  we highlight that, if $X_1$ exists in the set of the arguments of $p(X)$, both $NCQTE$ and $SCQTE$ can be useful as they can be robust against misspecification. Further, owing to the dimension-reduction structure, $SCQTE$ is worth of recommendation in large-dimensional scenarios. When $X_1$ is not in the set of the arguments of $p(X)$, all the estimators perform similarly in most cases.

\begin{table}[H]

  \centering
  \caption{The simulation results of $\hat \Delta_{\tau}(x_{1})$ under different scenarios}\label{tab1}
    \renewcommand\arraystretch{0.5}
  \scalebox{0.65}{
    \begin{tabular}{ccrrrrrrrrrrrrrrr}

    \toprule
    \multicolumn{17}{c}{n=500,    $Group\ 1: \{a=0.5,a_1=1.1,b_1=1.2\}$} \\
    \midrule
          &       & \multicolumn{3}{c}{OCQTE} &       & \multicolumn{3}{c}{PCQTE} &       & \multicolumn{3}{c}{NCQTE} &       & \multicolumn{3}{c}{SCQTE} \\
\cmidrule{3-5}\cmidrule{7-9}\cmidrule{11-13}\cmidrule{15-17}    Model & $x_1$ & \multicolumn{1}{c}{-0.2} & \multicolumn{1}{c}{0} & \multicolumn{1}{c}{0.2} &       & \multicolumn{1}{c}{-0.2} & \multicolumn{1}{c}{0} & \multicolumn{1}{c}{0.2} &       & \multicolumn{1}{c}{-0.2} & \multicolumn{1}{c}{0} & \multicolumn{1}{c}{0.2} &       & \multicolumn{1}{c}{-0.2} & \multicolumn{1}{c}{0} & \multicolumn{1}{c}{0.2} \\
\cmidrule{1-5}\cmidrule{7-9}\cmidrule{11-13}\cmidrule{15-17}    $Model~1$ & Bias  & 0.0213  & 0.0235  & 0.0209  &       & 0.0212  & 0.0233  & 0.0209  &       & 0.0305  & 0.0311  & 0.0269  &       & 0.0284  & 0.0298  & 0.0263  \\
          & SD    & 0.0343  & 0.0261  & 0.0315  &       & 0.0330  & 0.0249  & 0.0306  &       & 0.0325  & 0.0246  & 0.0304  &       & 0.0328  & 0.0249  & 0.0306  \\
          & MSE   & 0.0016  & 0.0012  & 0.0014  &       & 0.0015  & 0.0012  & 0.0014  &       & 0.0020  & 0.0016  & 0.0016  &       & 0.0019  & 0.0015  & 0.0016  \\
\cmidrule{1-5}\cmidrule{7-9}\cmidrule{11-13}\cmidrule{15-17}    $Model~2$ & Bias  & 0.0163  & 0.0183  & 0.0165  &       & 0.0162  & 0.0179  & 0.0164  &       & 0.0198  & 0.0217  & 0.0201  &       & 0.0243  & 0.0264  & 0.0244  \\
          & SD    & 0.0334  & 0.0259  & 0.0323  &       & 0.0329  & 0.0251  & 0.0316  &       & 0.0327  & 0.0248  & 0.0313  &       & 0.0326  & 0.0250  & 0.0314  \\
          & MSE   & 0.0014  & 0.0010  & 0.0013  &       & 0.0013  & 0.0009  & 0.0013  &       & 0.0015  & 0.0011  & 0.0014  &       & 0.0017  & 0.0013  & 0.0016  \\
\cmidrule{1-5}\cmidrule{7-9}\cmidrule{11-13}\cmidrule{15-17}    $Model~3$ & Bias  & 0.0613  & 0.0802  & 0.0829  &       & 0.0610  & 0.0805  & 0.0833  &       & 0.0812  & 0.1007  & 0.0984  &       & 0.0819  & 0.1024  & 0.1004  \\
          & SD    & 0.0432  & 0.0403  & 0.0439  &       & 0.0417  & 0.0378  & 0.0424  &       & 0.0411  & 0.0379  & 0.0425  &       & 0.0413  & 0.0375  & 0.0424  \\
          & MSE   & 0.0056  & 0.0081  & 0.0088  &       & 0.0055  & 0.0079  & 0.0087  &       & 0.0083  & 0.0116  & 0.0115  &       & 0.0084  & 0.0119  & 0.0119  \\
    \midrule
    \multicolumn{17}{c}{n=1000,    $Group\ 1: \{a=0.5,a_1=1.1,b_1=1.2\}$} \\
    \midrule
          &       & \multicolumn{3}{c}{OCQTE} &       & \multicolumn{3}{c}{PCQTE} &       & \multicolumn{3}{c}{NCQTE} &       & \multicolumn{3}{c}{SCQTE} \\
\cmidrule{3-5}\cmidrule{7-9}\cmidrule{11-13}\cmidrule{15-17}    Model & $x_1$ & \multicolumn{1}{c}{-0.2} & \multicolumn{1}{c}{0} & \multicolumn{1}{c}{0.2} &       & \multicolumn{1}{c}{-0.2} & \multicolumn{1}{c}{0} & \multicolumn{1}{c}{0.2} &       & \multicolumn{1}{c}{-0.2} & \multicolumn{1}{c}{0} & \multicolumn{1}{c}{0.2} &       & \multicolumn{1}{c}{-0.2} & \multicolumn{1}{c}{0} & \multicolumn{1}{c}{0.2} \\
\cmidrule{1-5}\cmidrule{7-9}\cmidrule{11-13}\cmidrule{15-17}    $Model~1$ & Bias  & 0.0193  & 0.0198  & 0.0191  &       & 0.0192  & 0.0199  & 0.0192  &       & 0.0276  & 0.0268  & 0.0243  &       & 0.0249  & 0.0252  & 0.0234  \\
          & SD    & 0.0236  & 0.0183  & 0.0223  &       & 0.0227  & 0.0174  & 0.0218  &       & 0.0225  & 0.0172  & 0.0214  &       & 0.0235  & 0.0180  & 0.0218  \\
          & MSE   & 0.0009  & 0.0007  & 0.0009  &       & 0.0009  & 0.0007  & 0.0008  &       & 0.0013  & 0.0010  & 0.0010  &       & 0.0012  & 0.0010  & 0.0010  \\
\cmidrule{1-5}\cmidrule{7-9}\cmidrule{11-13}\cmidrule{15-17}    $Model~2$ & Bias  & 0.0134  & 0.0138  & 0.0138  &       & 0.0132  & 0.0137  & 0.0137  &       & 0.0154  & 0.0162  & 0.0161  &       & 0.0195  & 0.0207  & 0.0200  \\
          & SD    & 0.0234  & 0.0191  & 0.0235  &       & 0.0230  & 0.0185  & 0.0230  &       & 0.0239  & 0.0188  & 0.0231  &       & 0.0230  & 0.0187  & 0.0230  \\
          & MSE   & 0.0007  & 0.0006  & 0.0007  &       & 0.0007  & 0.0005  & 0.0007  &       & 0.0008  & 0.0006  & 0.0008  &       & 0.0009  & 0.0008  & 0.0009  \\
\cmidrule{1-5}\cmidrule{7-9}\cmidrule{11-13}\cmidrule{15-17}    $Model~3$ & Bias  & 0.0555  & 0.0722  & 0.0775  &       & 0.0555  & 0.0719  & 0.0774  &       & 0.0749  & 0.0906  & 0.0914  &       & 0.0746  & 0.0905  & 0.0920  \\
          & SD    & 0.0317  & 0.0285  & 0.0304  &       & 0.0306  & 0.0275  & 0.0299  &       & 0.0302  & 0.0272  & 0.0298  &       & 0.0303  & 0.0271  & 0.0296  \\
          & MSE   & 0.0041  & 0.0060  & 0.0069  &       & 0.0040  & 0.0059  & 0.0069  &       & 0.0065  & 0.0089  & 0.0092  &       & 0.0065  & 0.0089  & 0.0093  \\
    \midrule
    \multicolumn{17}{c}{n=500,    $Group\ 2: \{a=0.5,a_1=0.9,b_1=1.1\}$} \\
    \midrule
          &       & \multicolumn{3}{c}{OCQTE} &       & \multicolumn{3}{c}{PCQTE} &       & \multicolumn{3}{c}{NCQTE} &       & \multicolumn{3}{c}{SCQTE} \\
\cmidrule{3-5}\cmidrule{7-9}\cmidrule{11-13}\cmidrule{15-17}    Model & $x_1$ & \multicolumn{1}{c}{-0.2} & \multicolumn{1}{c}{0} & \multicolumn{1}{c}{0.2} &       & \multicolumn{1}{c}{-0.2} & \multicolumn{1}{c}{0} & \multicolumn{1}{c}{0.2} &       & \multicolumn{1}{c}{-0.2} & \multicolumn{1}{c}{0} & \multicolumn{1}{c}{0.2} &       & \multicolumn{1}{c}{-0.2} & \multicolumn{1}{c}{0} & \multicolumn{1}{c}{0.2} \\
\cmidrule{1-5}\cmidrule{7-9}\cmidrule{11-13}\cmidrule{15-17}    $Model~1$ & Bias  & 0.0213  & 0.0235  & 0.0209  &       & 0.0212  & 0.0233  & 0.0209  &       & 0.0289  & 0.0295  & 0.0253  &       & 0.0280  & 0.0295  & 0.0260  \\
          & SD    & 0.0343  & 0.0261  & 0.0315  &       & 0.0330  & 0.0249  & 0.0306  &       & 0.0323  & 0.0242  & 0.0303  &       & 0.0328  & 0.0248  & 0.0306  \\
          & MSE   & 0.0016  & 0.0012  & 0.0014  &       & 0.0015  & 0.0012  & 0.0014  &       & 0.0019  & 0.0015  & 0.0016  &       & 0.0019  & 0.0015  & 0.0016  \\
\cmidrule{1-5}\cmidrule{7-9}\cmidrule{11-13}\cmidrule{15-17}    $Model~2$ & Bias  & 0.0163  & 0.0183  & 0.0165  &       & 0.0162  & 0.0179  & 0.0164  &       & 0.0186  & 0.0203  & 0.0187  &       & 0.0239  & 0.0261  & 0.0240  \\
          & SD    & 0.0334  & 0.0259  & 0.0323  &       & 0.0329  & 0.0251  & 0.0316  &       & 0.0327  & 0.0253  & 0.0317  &       & 0.0327  & 0.0250  & 0.0314  \\
          & MSE   & 0.0014  & 0.0010  & 0.0013  &       & 0.0013  & 0.0009  & 0.0013  &       & 0.0014  & 0.0011  & 0.0014  &       & 0.0016  & 0.0013  & 0.0016  \\
\cmidrule{1-5}\cmidrule{7-9}\cmidrule{11-13}\cmidrule{15-17}    $Model~3$ & Bias  & 0.0613  & 0.0802  & 0.0829  &       & 0.0610  & 0.0805  & 0.0833  &       & 0.0777  & 0.0968  & 0.0953  &       & 0.0806  & 0.1008  & 0.0992  \\
          & SD    & 0.0432  & 0.0403  & 0.0439  &       & 0.0417  & 0.0378  & 0.0424  &       & 0.0409  & 0.0375  & 0.0421  &       & 0.0413  & 0.0376  & 0.0423  \\
          & MSE   & 0.0056  & 0.0081  & 0.0088  &       & 0.0055  & 0.0079  & 0.0087  &       & 0.0077  & 0.0108  & 0.0109  &       & 0.0082  & 0.0116  & 0.0116  \\
    \midrule
    \multicolumn{17}{c}{n=1000,    $Group\ 2: \{a=0.5,a_1=0.9,b_1=1.1\}$} \\
    \midrule
          &       & \multicolumn{3}{c}{OCQTE} &       & \multicolumn{3}{c}{PCQTE} &       & \multicolumn{3}{c}{NCQTE} &       & \multicolumn{3}{c}{SCQTE} \\
\cmidrule{3-5}\cmidrule{7-9}\cmidrule{11-13}\cmidrule{15-17}    Model & $x_1$ & \multicolumn{1}{c}{-0.2} & \multicolumn{1}{c}{0} & \multicolumn{1}{c}{0.2} &       & \multicolumn{1}{c}{-0.2} & \multicolumn{1}{c}{0} & \multicolumn{1}{c}{0.2} &       & \multicolumn{1}{c}{-0.2} & \multicolumn{1}{c}{0} & \multicolumn{1}{c}{0.2} &       & \multicolumn{1}{c}{-0.2} & \multicolumn{1}{c}{0} & \multicolumn{1}{c}{0.2} \\
\cmidrule{1-5}\cmidrule{7-9}\cmidrule{11-13}\cmidrule{15-17}    $Model~1$ & Bias  & 0.0193  & 0.0198  & 0.0191  &       & 0.0192  & 0.0199  & 0.0192  &       & 0.0260  & 0.0253  & 0.0229  &       & 0.0244  & 0.0249  & 0.0232  \\
          & SD    & 0.0236  & 0.0183  & 0.0223  &       & 0.0227  & 0.0174  & 0.0218  &       & 0.0223  & 0.0169  & 0.0212  &       & 0.0243  & 0.0184  & 0.0218  \\
          & MSE   & 0.0009  & 0.0007  & 0.0009  &       & 0.0009  & 0.0007  & 0.0008  &       & 0.0012  & 0.0009  & 0.0010  &       & 0.0012  & 0.0010  & 0.0010  \\
\cmidrule{1-5}\cmidrule{7-9}\cmidrule{11-13}\cmidrule{15-17}    $Model~2$ & Bias  & 0.0134  & 0.0138  & 0.0138  &       & 0.0132  & 0.0137  & 0.0137  &       & 0.0141  & 0.0146  & 0.0148  &       & 0.0190  & 0.0203  & 0.0197  \\
          & SD    & 0.0234  & 0.0191  & 0.0235  &       & 0.0230  & 0.0185  & 0.0230  &       & 0.0256  & 0.0212  & 0.0251  &       & 0.0233  & 0.0188  & 0.0230  \\
          & MSE   & 0.0007  & 0.0006  & 0.0007  &       & 0.0007  & 0.0005  & 0.0007  &       & 0.0009  & 0.0007  & 0.0008  &       & 0.0009  & 0.0008  & 0.0009  \\
\cmidrule{1-5}\cmidrule{7-9}\cmidrule{11-13}\cmidrule{15-17}    $Model~3$ & Bias  & 0.0555  & 0.0722  & 0.0775  &       & 0.0555  & 0.0719  & 0.0774  &       & 0.0711  & 0.0870  & 0.0885  &       & 0.0732  & 0.0891  & 0.0907  \\
          & SD    & 0.0317  & 0.0285  & 0.0304  &       & 0.0306  & 0.0275  & 0.0299  &       & 0.0299  & 0.0269  & 0.0295  &       & 0.0302  & 0.0271  & 0.0295  \\
          & MSE   & 0.0041  & 0.0060  & 0.0069  &       & 0.0040  & 0.0059  & 0.0069  &       & 0.0060  & 0.0083  & 0.0087  &       & 0.0063  & 0.0087  & 0.0091  \\
    \bottomrule
      \end{tabular}}%
\end{table}%

\begin{figure}[!htbp]
\vspace*{-5pt}
\centering
  \includegraphics[width=13cm,height=7.6cm]{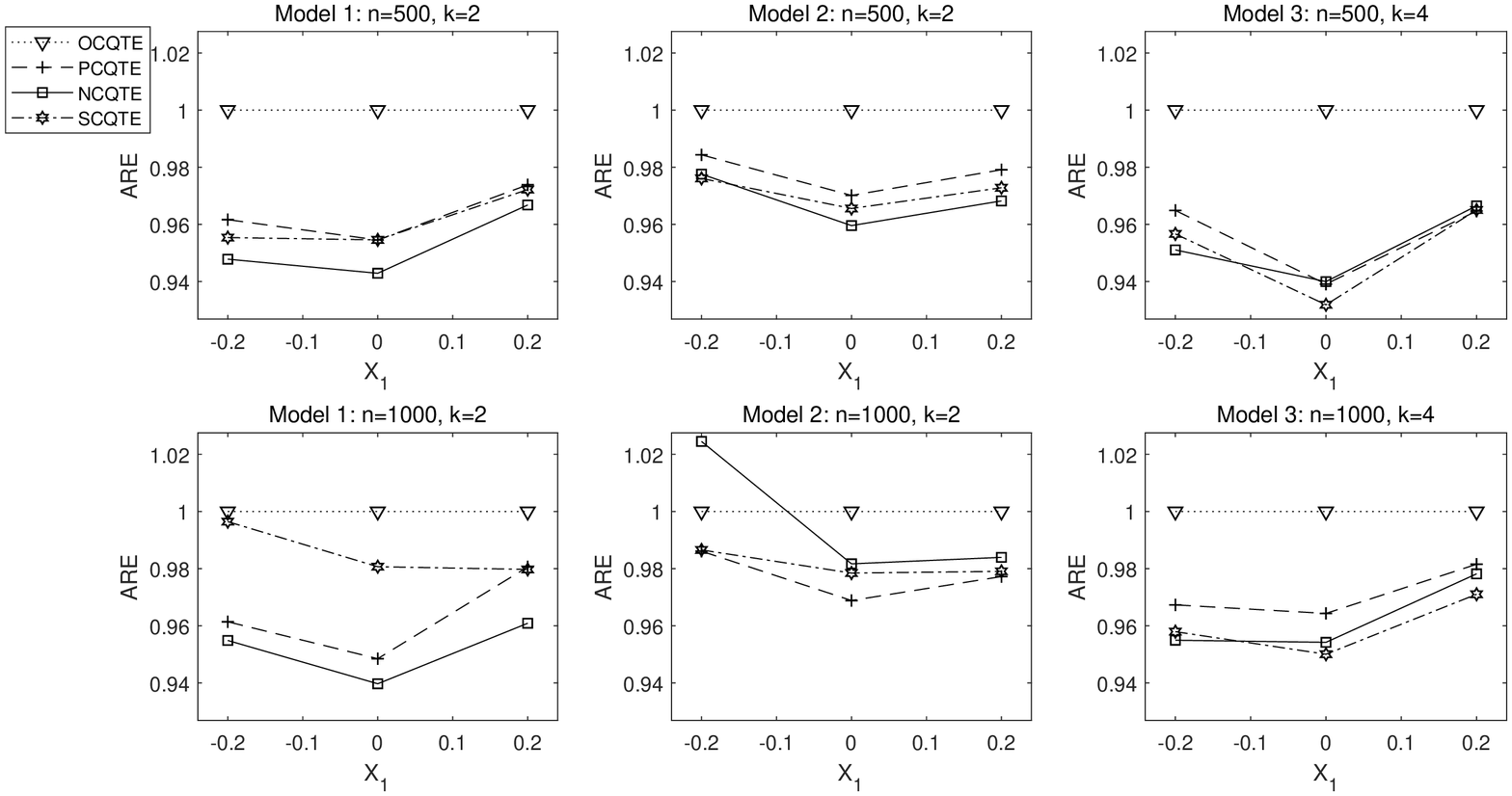}
  \vspace*{-10pt}
\caption{\scriptsize{ The asymptotic relative efficiency(ARE) about $SD$ against that of $OCQTE$ with the tuning parameters in Group~1.}}
 \label{group1}
\vspace*{-10pt}
\end{figure}
%
%
%
\section{A real data example}

In this section, we apply the proposed methodology to estimate the $CQTE$ function to investigate the quantile effect of maternal smoking on birth weight over the mother's age. We adopt a dataset based on the records between 1988 and
2002 by the North Carolina State Center Health Services, which can be obtained from Robert Lieli's website \url{http://www.personal.ceu.hu/staff/Robert_ Lieli/cate-birthdata.zip}. Noted that this dataset was also analyzed by \cite{abrevaya2015}, who aimed to estimate the conditional average treatment effect (CATE) of maternal smoking on birth weight by selecting $X_1$ as mother's age. Focusing on first-time mothers, they observed that the $CATE$ function is mostly negative, and \cite{abrevaya2015} also noted that $CATE$ is stronger( more negative) for older mothers.

We also choose mother's age as $X_1$, and aim to explore more information about the conditional smoking effect besides the average treatment effect provides. Both the low-birthweight (LBW) (weighing less than 2,500 grams) and high-birthweight(HBW) (weighing more than 4,000 grams) babies should receive attention in the literature. For example, we may also want to know, when the mother is older, whether the smoking effect will be stronger for the LBW babies, or whether there exist different trends of the smoking effect over mother's age for the LBW and HBW babies. Thus, we estimate the conditional quantile treatment effect $(CQTE)$, under $\tau=0.1,0.5,0.9$ respectively, to investigate how the quantile treatment effect varies with different values of  mother's age and different babies groups.

Before estimation, we first introduce some details and settings about the dataset. We restrict our sample to white and first-time mothers, thus the sample size is $n=433,558$ while the smoking sample size is $74,386$. The outcome $Y$ here is birth weight measured in grams and the treatment indicator variable $D$ is a binary variable. When $D=1$, it means the mother smokes and $D=0$ otherwise. Further, to ensure the unconfoundedness assumption, we choose a large set of variables as $X$, including the mother's age, education level, the month of the first prenatal visit (=10 if the prenatal care is foregone), the number of prenatal visits, and indicators for  baby's gender,  mother's marital status, whether or not the father's age is missing,  gestational diabetes,  hypertension, amniocentesis,  ultra sound exams, the previous (terminated) pregnancies, and alcohol use.

We estimate the $CQTE$ function $\Delta_\tau(x_1)$ in the interval between ages 15 and 35 under three different quantiles, i.e, $\Delta_{0.1}(x_1)$,$\Delta_{0.5}(x_1)$ and $\Delta_{0.9}(x_1)$ respectively. Since the dimension of $X$ is large, we use a semiparametric model for the propensity score that has a single index structure such that the dimensionality and model misspecification problems can be alleviated. Thus, we first use the sufficient dimension reduction method, $SIR$, to estimate the  index. However, in order to capture the nonlinear information of $p(x)$, the explanatory variables $X^{*}$ used in estimation consist of all the elements of $X$, the square of the mother's age, and the interaction terms between the mother's age and all other elements of $X$. When it comes to the selection of bandwidth, we set $h_2=\hat\sigma_dn^{-1/3}$ and $h=\hat \sigma_1n^{-1/5},$ where $\hat\sigma_d=\sqrt{var(\hat \alpha^{\top}x^*)}$, $\hat \alpha$ is the estimated linear index direction and $\hat\sigma_1=2\sqrt{var(x_1)}.$ As for kennel function, we use a regular Gaussian kernel as simulation studies.

Figure~\ref{scqte_fig1} displays the results of the estimated $CQTE(x_{1})$ as a function of the mother's age in the range of 15 to 35 years old. There are several points we want to highlight: (1) {\it The $CQTE(x_{1})$ for the effect of maternal smoking on birth weight is remarkably negative.} All three $CQTE(x_{1})$ curves range from about -140 grams to -300 grams, which means if a mother smokes during the pregnancy period, the birth weight of her baby will most likely decrease. This finding is in accordance with the conclusion of \cite{abrevaya2015}. (2) {\it The LBW babies suffer the most from maternal smoking and get thinner across the mother's age.} When we focus on $\hat \Delta_{0.1}(x_1)$ curve, it is at the bottom of all the three curves. Furthermore, we can observe it has a decreasing trend over mother's age. Thus this suggests that older mothers would be more urgently quit smoking to avoid the ultra-low-weight baby to occur. (3) {\it The trend of $\Delta_{0.5}(x_1)$ varies with the mother's age.} As for $\Delta_{0.5}(x_1)$, we can also find a decreasing trend from 16 to around 22 years of age, while the curve is rather stable between the age of 23 to 28. As the relationship between $median$ and $average$, it can be expected that $\Delta_{0.5}(x_1)$ is much like $CATE(x_1).$

%
\begin{figure}[!htbp]
\vspace*{-10pt}
\centering
\includegraphics[width=12cm,height=6.4cm]{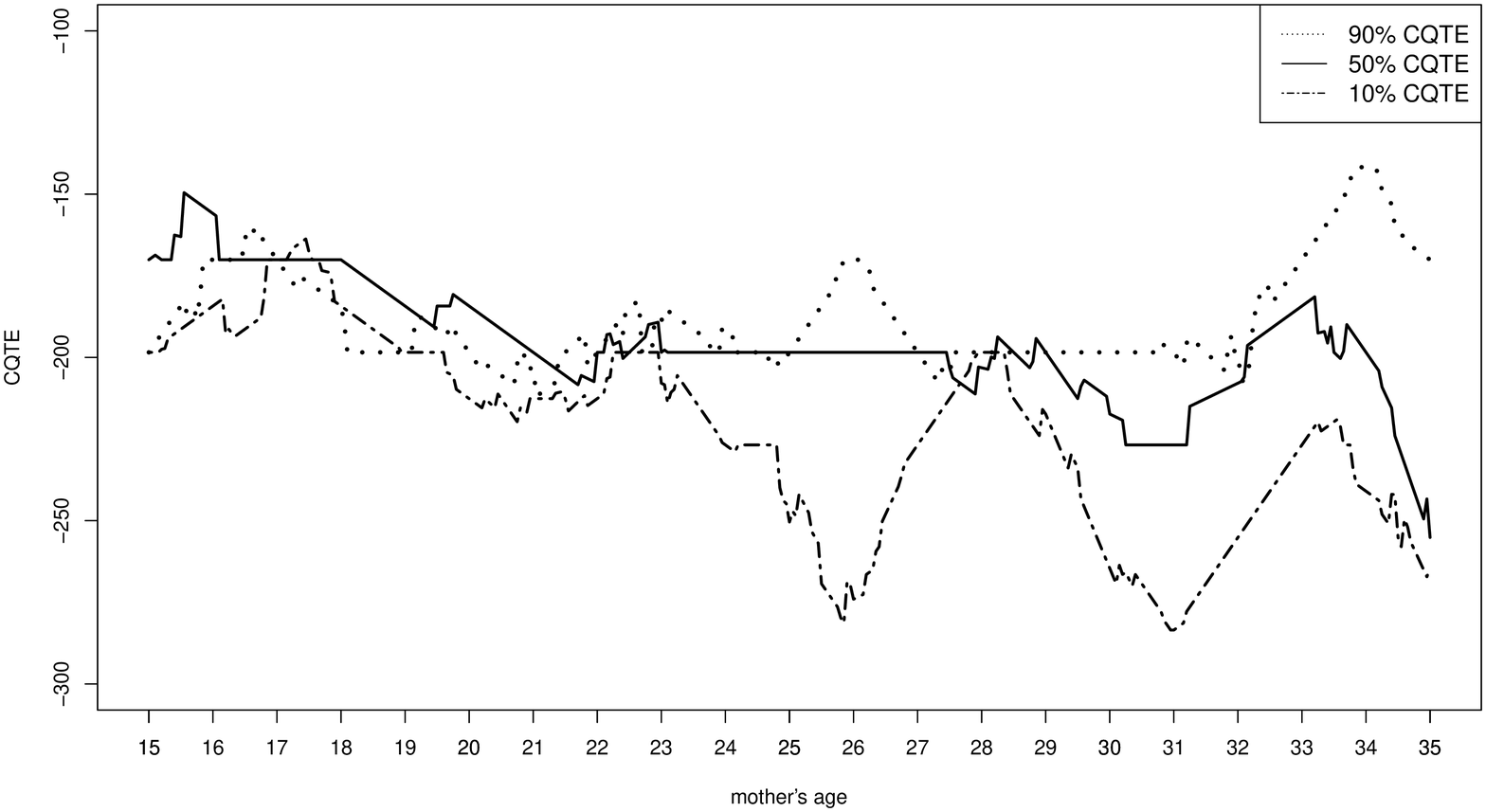}
\caption{\scriptsize{Three conditional quantile treatment effects (CQTE) curves over mother's age: $\hat\Delta_{0.9}(x_1)$ (dotted line), $\hat\Delta_{0.5}(x_1)$ (solid line) and $\hat\Delta_{0.1}(x_1)$ (dashed line).}}\label{scqte_fig1}
\vspace*{-8pt}
\end{figure}

%
%

\section{Conclusion}

In this paper, we propose the estimation of conditional
quantile treatment effect (CQTE), aimed to capture the conditional treatment effect in a specific subgroup. Four estimators are proposed when the propensity score is under  true function, parametric, nonparametric and semiparametric  structure: $OCQTE$, $PCQTE$, $NCQTE$ and $SCQTE$ where $OCQTE$ mainly serves as a benchmark for the comparison among the other three estimators. The asymptotic properties of the estimators are systematically investigated. The new findings show that the estimations under the conditional framework is rather different from their unconditional counterparts. More importantly,  under conditional framework, two factors play important role for the estimation efficiency of nonparametric and semiparametric-based estimators: 1) the convergence rate of the estimated propensity score; and 2) the affiliation  of the given covariates to the set of  arguments of the propensity score. These are not the cases for the unconditional counterparts in studying treatment effect. One more issue is about semiparametric efficiency the unconditional counterparts can achieve when nonparamatric estimation is used for propensity score. Under the conditional framework, $CQTE$s are  functions of the given covariates, it is unclear what would be defined as a semiparametric efficiency. It seems to involve uniformly asymptotic efficiency over a function. Thus, we leave it to a further study.

{Note that the two-step estimation procedure of $CQTE(X_1)$ could be extended to deal with  a more general treatment effect function:
  \begin{eqnarray}\label{mqte}
  \begin{split}
 &&M(X_1)=\arg\min_{a}
 E\left[\psi(Y(1),a)\mid
 X_{1}\right]-\arg\min_{a}
 E\left[\psi(Y(0),a)\mid
 X_{1}\right],\\
 &&=\arg\min_{a}
 E\left[\frac{D}{p(X)}\psi(Y,a)\mid
 X_{1}\right]-\arg\min_{a}
 E\left[\frac{1-D}{1-p(X)}\psi(Y,a)\mid
 X_{1}\right].
 \end{split}
 \end{eqnarray}
Here $\psi(\cdot)$ is a known real-value function. When the loss function $\psi(Y,a)$ equals  $\rho_\tau(Y-a)$ or $(Y-a)^2$, we can identify $CQTE(X_1)$ or $CATE(X_1)$  accordingly. The estimator $\hat M(X_1)$ could be obtained by solving the sample analogy  of (\ref{mqte}) and we could similarly derive the asymptotic behaviors of $\hat M(x_{10}).$ As it would involve  different optimization issues and theoretical investigation, we will give a detailed research in a later study.

Further,  in this field, potential outcome regression and doubly robust estimation are also the basic methodologies, the relevant studies are worthwhile.  However, as this paper mainly focuses on a systematic investigation on the asymptotic efficiencies of different propensity score-based estimations, the systematic studies  about potential outcome regression and doubly robust estimation will be the topics in the near future.
 }

%
%
%
\vskip 14pt
\noindent {\large\bf Supplementary Materials}
  This Supplementary Material contains with the technical lemma and
proofs of the main results.
\par

\bibhang=1.7pc
\bibsep=2pt
\fontsize{9}{14pt plus.8pt minus .6pt}\selectfont
\renewcommand\bibname{\large \bf References}
\expandafter\ifx\csname
natexlab\endcsname\relax\def\natexlab#1{#1}\fi
\expandafter\ifx\csname url\endcsname\relax
  \def\url#1{\texttt{#1}}\fi
\expandafter\ifx\csname urlprefix\endcsname\relax\def\urlprefix{URL}\fi

\lhead[\footnotesize\thepage\fancyplain{}\leftmark]{}\rhead[]{\fancyplain{}\rightmark\footnotesize{} }

\vskip .65cm
%
%
\end{document}


\def\spacingset#1{\renewcommand{\baselinestretch}%
{#1}\small\normalsize} \spacingset{1.5}
	\date{}

\if0\blind
{
  \title{\bf Supplement of ``The Role of Propensity Score Structure in Asymptotic Efficiency of Estimated Conditional Quantile Treatment
Effect''}
  \author{Niwen Zhou\\
    School of Statistics, Beijing Normal University\\
    and \\
    Xu Guo\\
   School of Statistics, Beijing Normal University\\
    and \\
    Lixing Zhu \\
   Department of Mathematics, Hong Kong Baptist University}
  \maketitle
} \fi

\if1\blind
{
  \bigskip
  \bigskip
  \bigskip
  \begin{center}
    {\LARGE\bf  Supplement of ``The Role of Propensity Score Structure in Asymptotic Efficiency of Estimated Conditional Quantile Treatment
Effect''}
\end{center}
  \medskip
} \fi

\spacingset{1.5} 

 \def\theequation{S.\arabic{section}.\arabic{equation}}
\def\thesection{S.\arabic{section}}

\section{Notation}

Without loss of
generality, we consider the situation $X_1 \in R$ in the following proofs for ease of exposition.

Before we present the proof, we first give some notations in this section.
\begin{itemize}

\item[(1)] $Y_i^*=Y_i- \tilde q_{1,\tau}(X_{1i}),$ $\tilde q_{1,\tau}(X_{1i})=q_{1,\tau}(x_{10})
-q'_{1,\tau}(x_{10})(X_{1i}-x_{10}),X_{1i}^*=(1,\frac{(X_{1i}-x_{10})}{h})^{\top}$;

\item[(2)] $\theta_{\tau}:=(u_\tau,v_\tau):=\sqrt{nh}
    \left(a-q_{1,\tau}(x_{10}),h(b-q'_{1,\tau}(x_{10}))\right)$;

\item[(3)] $\eta_{i}^*=\mathbb I\left(Y_i^*\leq
0\right)-\tau$, \, $m_{j,\tau}(X)
=\frac{E[\mathbb I(Y(j)\leq \tilde q_{j,\tau}(X_{1i}))-\tau\mid
X]}{f_j(q_{j,\tau}(x_{10})\mid x_{10})}$,$j=0,1;$

\item[(5)]  $\psi(p(X_i),Z_i)=
\frac{D_i}{p(X_i)} \eta_{1,
\tau}(Y_i) -\frac{1-D_i}{1-p(X_i)} \eta_{0,\tau}(Y_i)$,
$\eta_{j,\tau}(Y_i)=\frac{\mathbb I(Y^*_i\leq
0)-\tau}{f_j(q_{j,\tau}(x_{10})\mid x_{10})},j=0,1;$

\item[(6)] $K_i=K(\frac{X_{1i}
-x_{10}}{h})$, $\mu_{s_1}(K)=\int u^{s_1}K(u)du,$ $\parallel
K\parallel_2^2=\int K(u)^2 du$;
\item[(7)] $Z_i=(X_i,D_i,Y_i)$,
$f(X_1)$ is the density function of $X_1$, $f_x(X)$ is the density function of $X$, $f_\alpha(\alpha{\top}X)$ is the density function of $\alpha^{\top}X$, and denote $f_j(y(j)\mid x_{1})$ is the point $(Y(j)=y(j), X_1=x_{1})$ in
the conditional density function of $(Y(j)\mid X_1)$, $j=0,1;$
\item[(8)] $C$ stands for a generic bounded constant, $\Lambda$ is the $\sigma$-field generated by  $X_{11},\cdots,X_{1n}.$
\end{itemize}

\section{Proofs}
\counterwithout{equation}{section}

\subsection{Proof of Theorem~2.1: under the true $p(x)$.}

Since the arguments for $q_{1,\tau}(x_{10})$ and
$q_{0,\tau}(x_{10}) $ are similar, in the following we focus on the
estimation of $q_{1,\tau}(x_{10})$. Recall the objective function is:
\begin{eqnarray*}
G_{\tau}(a,b,p)&=&\sum_{i=1}^{n}\frac{D_i}{p(X_i)}
\rho_{\tau}(Y_i-a-b(X_{1i}-x_{10}))K_i.
\end{eqnarray*}
Obviously, minimizing $G_{\tau}(  a,b,p)$ is equivalent to
minimizing $A_{n,2}$ where
\begin{eqnarray*}
A_{n,2}= G_{\tau}(  a,b,p)-
G_{\tau}\left(q_{1,\tau}(x_{10}),q'_{1,\tau}(x_{10}),p\right).
\end{eqnarray*}
{It is easy to see that we cannot easily get a closed form for the minimizers of the above introduced objective function.
Thus to easily analyse the asymptotic behaviors of the estimators, we follow the arguments in the proof of Theorem~1 in \cite{carroll1997} to first obtain a quadratic approximation to the objective function $A_{n,2}$. The approximation error is $o_p(1)$. Then using the convexity lemma (Pollard 1991), the minimizers of the objective function, $\theta_{\tau}=\sqrt{nh}\Big(a-q_{1,\tau}(x_{10}),h\big(b-q'_{1,\tau}(x_{10})\big)\Big)$, can also be well approximated by the minimizers of the quadratic approximation. This approach has been adopted by many authors, see for instance \cite{firpo2007} and \cite{kai2010}. This approach is effective, especially when the objective function is not sufficiently smooth.}
We  now give the detail. Rewrite $A_{n,2}$ as:
\begin{eqnarray*}
A_{n,2}=\sum_{i=1}^{n}\frac{D_i}{ {p}(X_i)}
\left[\rho_{\tau}(-\theta_\tau^{\top}X_{1i}^*/\sqrt{nh}
+Y_i^*)-\rho_{\tau}(Y_i^*)\right]K_i.
\end{eqnarray*}
Applying the identity \cite{knight1998}:
\begin{equation*}
\rho_{\tau}(x-y)-\rho_{\tau}(x)=y{\mathbb \{I(x\leq
0)}-\tau\}+
\int_0^y \mathbb I\{x\leq z\}-\mathbb I\{x\leq 0\}dz,
\end{equation*}
we can transform $A_{n,2}$ to be
\begin{eqnarray*}
\begin{split}
A_{n,2}&=
\frac{1}{\sqrt{nh}}\sum_{i=1}^{n}\frac{D_i\theta_\tau^{\top}X_{1i}^*
}{ {p}(X_i)}K_i\left(\mathbb I (Y_i^*\leq 0)-\tau \right) +\sum_{i=1}^{n}\frac{D_i}{
{p}(X_i)}K_i\int_0^{\theta_\tau^{\top}X_{1i}^*/\sqrt{nh}}\left[
\mathbb I(Y_i^*\leq t )-\mathbb I\left(Y_i^*\leq
0 \right) \right]dt.
\end{split}
\end{eqnarray*}
Define $R_{n,i}(Y_i,\theta_\tau)$ as
\begin{eqnarray*}
R_{n,i}(Y_i,\theta_\tau)=\int_0^{\theta_\tau^{\top}X_{1i}^*/\sqrt{nh}}\left[
\mathbb I(Y_i^*\leq t )-\mathbb I\left(Y_i^*\leq
0 \right) \right]dt.
\end{eqnarray*}
Next, we get an approximation of $A_{n,2}$ as
\begin{eqnarray}
A_{n,2}
=\theta_\tau^{\top}W_n+\sum_{i=1}^{n}\frac{D_i}{ {p}(X_i)}K_iR_{n,i}(Y_i,\theta_\tau),
\end{eqnarray}
where $W_n=\frac{1}{\sqrt{nh}}\sum_{i=1}^{n}
\frac{D_i}{p(X_i)}X_{1i}^* K_i (\mathbb I\left(Y_i^*\leq
0 \right)-\tau).$
Some elementary calculations yield
\begin{eqnarray}
\frac{1}{nh}\sum_{i=1}^{n}\frac{D_i}{ {p}(X_i)}K_iR_{n,i}(Y_i,\theta_\tau)=\frac{1}{2}\theta_\tau ^{\top} S\theta_\tau+o_p(1),
\end{eqnarray}
 where $S=:f_1( q_{1,\tau}(x_{10})\mid x_{10})f(x_{10})\begin{pmatrix}
1&0\\
0&\mu_2(K)
\end{pmatrix}$.
It follows that
\begin{eqnarray}{\label{ocqte_case}}
\begin{split}
A_{n,2}
&=\theta_\tau^{\top}W_n+\frac{1}{2}\theta_\tau ^{\top}S\theta_\tau+o_p(1),\\
&=u_\tau W_{n1}+\frac{1}{2}u_\tau^2f_1( q_{1,\tau}(x_{10})\mid x_{10})f(x_{10})+v_\tau^{\top}W_{n2}\\
&\quad +
\frac{1}{2}f_1( q_{1,\tau}(x_{10})\mid x_{10 })f(x_{10})v_\tau^{\top}\mu_2(K)v_\tau+o_p(1).
\end{split}
\end{eqnarray}

Further, the application of  the convexity lemma \citep{pollard1991} leads to
\begin{eqnarray*}
& &\sqrt{nh}\left(\hat q^{ocqte}_{1,\tau }(x_{10})- {q_{1,\tau}
(x_{10})}\right)=-\frac{1}{\sqrt{nh}}
\frac{1}{f(x_{10})}\sum_{i=1}^{n}
\frac{D_i}{p(X_i)} \eta_{1,
\tau}(Y_i)K_i+o_p(1),
\end{eqnarray*}
Similarly, we can prove that
\begin{eqnarray*}
 \sqrt{nh}\left(\hat q^{ocqte}_{0,\tau }(x_{10})- {q_{0,\tau}
(x_{10})}\right)=-\frac{1}{\sqrt{nh}}
\frac{1}{f(x_{10})}\sum_{i=1}^{n}
\frac{1-D_i}{1-p(X_i)} \eta_{0,
\tau}(Y_i)K_i+o_p(1).
\end{eqnarray*}

Note that the definition of $CQTE$: $\Delta_\tau(x_{10})=q_{1,\tau}(x_{10})-q_{0,\tau}(x_{10})$. It can easily get that
\begin{eqnarray*}
& &\sqrt{nh}\left(\hat \Delta^{ocqte}_{\tau }(x_{10})-
{\Delta_{\tau}
(x_{10})}\right)=-\frac{1}{\sqrt{nh^l}}\frac{1}{f(x_{10})}
\sum_{i=1}^{n}\psi(p(X_i),Z_i)K_i+o_p(1),
\end{eqnarray*}

Next we   study the asymptotic bias and variance
of $\hat \Delta^{ocqte}_\tau(x_{10}).$

(1) {\it Bias of $\hat \Delta^{ocqte}_\tau(x_{10})$. }

Note that $X=(X_1,X_{(2)})$, and $\tau=E(\mathbb I(Y_i(1)\leq q_{1,\tau}(X_{1i}))\mid X_{1i})=F_1(q_{1,\tau}(X_{1i})\mid X_{1i})$. We have
\begin{eqnarray*}
\begin{split}
&Bias(\hat \Delta^{ocqte}_\tau(x_{10})\mid \Lambda)=-\frac{1}{{nh}}
\frac{1}{f(x_{10})}\sum_{i=1}^{n}
E\left(\phi_1(p(X_i),Z_i) \mid X_{1i}\right)K_i\\
&=-\frac{1}{{nh}}
\frac{1}{f(x_{10})}\sum_{i=1}^{n}\bigg(
\frac{F_1(\tilde q_{1,\tau}(X_{1i})\mid X_{1i})-F_1(q_{1,\tau}(X_{1i})\mid X_{1i})}
{f_1(q_{1,\tau}(x_{10})\mid x_{10})} \\
&\qquad\qquad\qquad\qquad\qquad-\frac{F_0(\tilde q_{0,\tau}(X_{1i})\mid X_{1i})-F_0(q_{0,\tau}(X_{1i})\mid X_{1i})}
{f_0(q_{0,\tau}(x_{10})\mid x_{10})}\bigg)K_i
\end{split}
\end{eqnarray*}
Recall that $\tilde q_{1,\tau}(X_{1i})=q_{1,\tau}(x_{10})+q'_{1,\tau}(x_{10})(X_{1i}-x_{10}) $ .
  By a Taylor expansion around $\tilde q_{1,\tau}(X_{1i})$, it follows that, if $K$ is a kernel with order $s_1=2$,
\begin{eqnarray*}
\begin{split}
Bias(\hat \Delta^{ocqte}_\tau(x_{10})\mid \Lambda)
&=\frac{q^{''}_{1,\tau}(x_{10})}{{2}}
\frac{1}{f(x_{10})f_1(q_{1,\tau}(x_{10})\mid x_{10})}\frac{1}{nh}\sum_{i=1}^{n}
{f_1(\tilde q_{1,\tau}(X_{1i})\mid X_{1i})(X_{1i}-x_{10})^2}
{}K_i \\
&\qquad -\frac{q''_{0,\tau}(x_{10})}{{2}}
\frac{1}{f(x_{10})f_0(q_{0,\tau}(x_{10})\mid x_{10})}\frac{1}{nh}\sum_{i=1}^{n}
{f_0(\tilde q_{0,\tau}(X_{1i})\mid X_{1i})(X_{1i}-x_{10})^2}
{}K_i \\
&=\frac{1}{2}\Delta_\tau''(x_{10})
\mu_2(K)h^2+o_p(h^2),
\end{split}
\end{eqnarray*}
where $\mu_2(K)=\int u^2K(u)du.$ While $s_1>2$, it can similarly get that $Bias(\hat \Delta^{ocqte}_\tau(x_{10})\mid \Lambda)=O_p(\mu_{s_1}h^{s_1}).$

(2) {\it Variance of $\hat \Delta^{ocqte}_\tau(x_{10})$. }

Let $\sigma_{ocqte}^2(X_1)=E(\psi^2(p(X),Z)\mid X_1).$ We have
\begin{eqnarray*}
\begin{split}
&Var(\hat\Delta^{ocqte}_{\tau}(x_{10}))=\frac{1}{nh^2}
\frac{1}{f^2(x_{10})}\big
\{E\left[\psi^2
K^2\left(\frac{X_1-x_{10}}{h}\right)  \right] -\left[E \left(\psi
K\left(\frac{X_1-x_{10}}{h}\right)\right)\right]^2
     \bigg\}\\
&=\frac{1}{nh^2}
\frac{1}{f^2(x_{10})}\left\{E\left[\sigma_{ocqte}^2(X_1)K^2\left(\frac{X_1-x_{10}}{h}\right)
\right]
\right\}+o\left(\frac{1}{nh}\right)\\
&=\frac{1}{nh}\frac{\parallel
K\parallel_2^2\sigma_{ocqte}^2(x_{10})}{f(x_{10})}+o\left(\frac{1}{nh}\right),
\end{split}
\end{eqnarray*}
where $\parallel K\parallel_2^2=\int K^2(t)dt.$
Under the assumptions
\begin{eqnarray}
\begin{split}
&\sigma_{ocqte}^2(x_{10})=E(\phi_1^2(p(X),Z)\mid X_1=x_{10})=E\bigg(\frac{E\big[(\mathbb I(Y(1)\leq \tilde q_{1,\tau}(X_{1}))-\tau)^2\mid X\big]}{p(X)f^2_1(q_{1,\tau}(x_{10})\mid x_{10})}\nonumber\\
&\qquad \qquad\qquad\quad\qquad\qquad \qquad \qquad \qquad+\frac{E\big[(\mathbb I(Y(0)\leq \tilde q_{0,\tau}(X_{1}))-\tau)^2\mid X\big]}{(1-p(X))f^2_0(q_{0,\tau}(x_{10})\mid x_{10})} \mid X_1=x_{10}\bigg).
\end{split}
\end{eqnarray}
 Under the Assumptions 1 through 5, we can get
 \begin{eqnarray*}
\sqrt{nh}\left(\hat {\Delta}^{ocqte}_\tau (x_{10})-\Delta_\tau
(x_{10})-b_{0}(x_{10})\right)
\stackrel{\mathcal{D}}{\longrightarrow}N\left(0,\frac{\parallel
K \parallel_2^2\sigma_{ocqte}^2(x_{10})}
{f(x_{10})}\right)
\end{eqnarray*}
where $b_0(x_{10})=O_p(\mu_{s_1}h^{s_1}).$
Thus we complete the proof of $Theorem$~2.1.
 \hfill$\Box$

\subsection{Proofs of The  Theorems and Corollaries under estimated $\hat p(x)$.}
As the proofs for these theorems and corollaries have some overlaps technically, we then give the preliminary of the proofs such that the proofs for the individual theorems can share.

\subsubsection{Preliminary of the proofs.}
To make the proofs smoothly, we first decompose the objective function as follows.
\begin{eqnarray*}
&&G_{\tau}( a,b,\hat p)\\
&=&\sum_{i=1}^{n}\left(\frac{D_i}{ {\hat
p}(X_i)}-\frac{D_i}{ {
p}(X_i)}\right)\left[\rho_{\tau}(Y_i-a-b(X_{1i}-x_{10}))-
\rho_{\tau}(Y_i-\tilde q_{1,\tau}(X_{1i}))\right]K(\frac{X_{1i}-x_{10}}{h})\nonumber\\
&&+\sum_{i=1}^n\frac{D_i}{{
p}(X_i)}\left[\rho_{\tau}(Y_i-a-b(X_{1i}-x_{10}))
-\rho_{\tau}(Y_i-\tilde q_{1,\tau}(X_{1i}))\right]K(\frac{X_{1i}-x_{10}}{h})\nonumber\\
&&+\sum_{i=1}^n\frac{D_i}{{ \hat
p}(X_i)}\rho_{\tau}(Y_i-\tilde q_{1,\tau}(X_{1i}))K(\frac{X_{1i}-x_{10}}{h})\nonumber\\
&:=&A_{n,1}+A_{n,2}+A_{n,3}.
\end{eqnarray*}

When it comes to $A_{n,2}$ and $A_{n,3}$, since $A_{n,2}$ is based on true $p(x)$, which has been discussed in Lemma~2.1, it follows immediately that
\begin{eqnarray*}
&A_{n,2}=u_\tau \frac{1}{\sqrt{nh}}\sum_{i=1}^{n}\frac{D_i}{p(X_i)}
(\mathbb I(Y_i\leq
\tilde q_{1,\tau}(X_{1i}))-\tau)K_i+\frac{1}{2}u_\tau^2f_1( q_{1,\tau}(x_{10})\mid x_{10})f(x_{10})\\
&+v_\tau^{\top}W_{n2}+
\frac{1}{2}f_1( q_{1,\tau}(x_{10})\mid x_{10})f(x_{10})v_\tau^{\top}\mu_2(K)v_\tau+o_p(1).
\end{eqnarray*}
Recall that in the objective function $G_{\tau}( a,b,\hat p)$, $A_{n,3}$ is free of the parameter of interest and thus the minimization does not involve this term. Thus, the following proofs focus on $A_{n,1}$.

Applying the identity \citep{knight1998}
\begin{equation*}
\rho_{\tau}(x-y)-\rho_{\tau}(x)=y{\mathbb \{I(x\leq
0)}-\tau\}+
\int_0^y \mathbb I\{x\leq z\}-\mathbb I\{x\leq 0\}dz,
\end{equation*}
we have
\begin{align}\label{rn}
&&\quad\quad\qquad\rho_{\tau}(Y_i-a-b(X_{1i}-x_{10}))-\rho_{\tau}(Y_i
-\tilde q_{1,\tau}(X_{1i})
)=\rho_{\tau}(-\theta_\tau^{\top}X_{1i}^*/\sqrt{nh}
+Y_i^*)-\rho_{\tau}(Y_i^*)\nonumber\\
&&=\frac{\theta_\tau^{\top}X_{1i}^*}{\sqrt{nh}}{\mathbb \{I(Y_i^*\leq
0)}-\tau\}+
\int_0^{\frac{\theta_\tau^{\top}X_{1i}^*}{\sqrt{nh}}} \big(\mathbb I\{Y_i^*\leq
t\}-\mathbb I\{Y_i^*\leq
0\}\big)dt \nonumber\qquad\qquad\qquad\qquad\\
&&:=\frac{\theta_\tau^{\top}X_{1i}^*}
{\sqrt{nh}}\eta_i^*+R_{n,i}(Y_i,\theta_\tau)\qquad\qquad\qquad\qquad\qquad\qquad\qquad\qquad\qquad\qquad\qquad\qquad.
\end{align}
Thus
\begin{eqnarray}
A_{n,1}&=&\sum_{i=1}^{n}\left(\frac{D_i}{ {\hat
p}(X_i)}-\frac{D_i}{ {
p}(X_i)}\right)[\frac{\theta_\tau^{\top}X_{1i}^*}{\sqrt{nh}}\eta_i^*+R_{n,i}(Y_i,\theta_\tau)]
K_i\nonumber\\
&=&\frac{\theta_\tau^{\top}}{\sqrt{nh}}
\sum_{i=1}^n[\phi(Z_i,\hat
p(X_i))-\phi(Z_i,p(X_i))]K_i+C_{n,2}:=\theta_\tau^{\top}
C_{n,1}+C_{n,2}.
\end{eqnarray}
Here $\phi(p(X_i),Z_i)=\frac{D_iX_{1i}^*}{p(X_i)}(\mathbb I\left(Y_i^*\leq
0\right)-\tau),C_{n,1}=\frac{1}{\sqrt{nh}}
\sum_{i=1}^n[\phi(Z_i,\hat p(X_i))-\phi(Z_i,p(X_i))]K_i,$
and
$C_{n,2}=\sum_{i=1}^n\left(\frac{D_i}{ {\hat
p}(x_i)}-\frac{D_i}{ {
p}(x_i)}\right)R_{n,i}(Y_i,\theta_\tau)K_i.$

 Further, if $G_{\tau}( a,b,\hat p)$ can be written as a function $\tilde G_{\tau}( a,b, p)+o_p(1)$ the solution of $G_{\tau}( a,b,\hat p)$ can have the same limiting properties of $\tilde G_{\tau}( a,b, p).$ See e.g.  Kai, Li and Zou (2010). In the following, we  prove Theorems~2.2 through 2.4
via first showing that $C_{n,2}=o_p(1)$ and discovering the asymptotic
expansion for $C_{n,1}$. Afterwards, we derive the asymptotic results about the solutions. Since the proofs are based on
three different propensity score estimators, thus we discuss
the proofs separately.

\subsubsection{Proof of Theorem~2.2.}
Under  Assumption~6, $\hat p(X_i)$ is a parametric estimator such that
$$\mathop{\sup}_{x\in\chi}\mid \hat p(X) -p(X)  \mid=O_p(n^{-1/2}).$$
Now we have:
\begin{eqnarray}
|C_{n,1}|&=&|
\frac 1{\sqrt{nh}}\sum_{i=1}^n\phi_p(\tilde p(X_i),Z_i)K_i(\hat
p(X_i)-p(X_i))|\nonumber\\
&\leq& \sqrt{nh} \mathop{\sup}_{x}\mid \hat p(X) -p(X)\mid
\frac{1}{nh}\sum_{i=1}^{n}\mid \phi_p(\tilde p(X_i),Z_i)K(\frac{X_{1i}-x_{10}}{h})\mid
\end{eqnarray}
Here $\tilde p(X_i)$ is between $p(X_i)$ and $\hat p(X_i)$.
Note that $h\rightarrow0$ and $\frac{1}{nh}\sum_{i=1}^{n}\mid
\phi_p(\tilde
p(X_i),Z_i)K(\frac{X_{1i}-x_{10}}{h})\mid=O_p(1).$ Thus we have
$C_{n,1}=O_p(\sqrt h)=o_p(1).$

Next, we  prove that $C_{n,2}=o_p(1).$
Clearly, $p(X)$ and $K_i$ are bounded, and thus
\begin{eqnarray}\label{theo2}
\begin{split}
|C_{n,2}|&=|\sum_{i=1}^n
\left(\frac{D_i}{ {\hat p}(X_i)}-\frac{D_i}{ {
p}(X_i)}\right)(K_iR_{n,i}(Y_i,\theta_\tau))\\
&\leq C  n\mathop{\sup}_{x}\mid \hat p(X) -p(X) \mid
\frac 1n\sum_{i=1}^n \mid K_iR_{n,i}(Y_i,\theta_\tau)\mid \\
\end{split}
\end{eqnarray}
where $C$ is a bounded constant, whose value may change for each appearance. By Chebyshev's Inequality, for any $c>0$ we bound the following probability as
\begin{eqnarray*} P(\frac 1n\sum_{i=1}^n \mid K_iR_{n,i}(Y_i, \theta_\tau)\mid>c)
&\le&  E(\mid
K_1R_{n,1}(Y_1,\theta_\tau)\mid)/c\\
&\le & \sqrt {E(\mid
K_1R_{n,1}(Y_1,\theta_\tau)\mid)^2}/c.
\end{eqnarray*}
Recall the definition of $R_{n,1}(Y_1,\theta_\tau)$ in (\ref{rn}), we can easily derive that $E(\mid
K_1R_{n,1}(Y_1,\theta_\tau)\mid)^2=O(n^{-3/2}h^{-1/2}).$ See a relevant reference  \cite{firpo2007}. Together with (\ref{theo2}), we have that
\begin{eqnarray}\label{cn2}|C_{n, 2}|=O_p((n/h)^{1/4}n^{-1/2})=o_p(1).
\end{eqnarray}
The above results imply that $$G_\tau(a,b,\hat
p)=A_{n,2}+A_{n,3}+o_p(1).$$
Since $A_{n,3}$ does not depend on $a$, we can get the
asymptotically linear representation  similarly by the proof
of Lemma~2.1. In other words,
\begin{eqnarray*}
& &\sqrt{nh}\left(\hat q^{pcqte}_{1,\tau }(x_{10})- {q_{1,\tau}
(x_{10})}\right)=-\frac{1}{\sqrt{nh}}
\frac{1}{f(x_{10})}\sum_{i=1}^{n}
\frac{D_i}{p(X_i)} \eta_{1,
\tau}(Y_i)K_i+o_p(1),
\end{eqnarray*}
By an analogous proof for $\hat q^{pcqte}_{1,\tau}(x_{10})$, we can get the asymptotical linear representation of $\hat
q^{pcqte}_{0,\tau}(x_{10})$ as
\begin{eqnarray*}
& &\sqrt{nh}\left(\hat q^{pcqte}_{0,\tau }(x_{10})- {q_{0,\tau}
(x_{10})}\right)=-\frac{1}{\sqrt{nh}}
\frac{1}{f(x_{10})}\sum_{i=1}^{n}
\frac{1-D_i}{1-p(X_i)} \eta_{0,
\tau}(Y_i)K_i+o_p(1),
\end{eqnarray*}
Recall that $\hat
\Delta_\tau(x_{10})=q_{1,\tau}(x_{10})-q_{0,\tau}(x_{10}).$
Thus the asymptotically linear representation of
$\hat \Delta^{pcqte}_\tau(x_{10})$ is
\begin{eqnarray*}
& &\sqrt{nh}\left(\hat \Delta^{pcqte}_{\tau }(x_{10})-
{\Delta}_{\tau} (x_{10})\right)=-\frac{1}{\sqrt{nh}}
\frac{1}{f(x_{10})}\sum_{i=1}^{n}\phi_1(p(X_i),Z_i)K_i+o_p(1),\\
\end{eqnarray*}
where $ \phi_1(p(X_i),Z_i)=\psi(p(X_i),Z_i)=
\frac{D_i}{p(X_i)} \eta_{1,
\tau}(Y_i) -\frac{1-D_i}{1-p(X_i)} \eta_{0,\tau}(Y_i).$

Obviously, $PCQTE$ has the same asymptotically linear representation of $OCQTE$, thus the rest proof of $Theorem$~2.2 follows from that of $Theorem$~2.1 immediately.
Thus we can complete the proof of $Theorem$~2.2.
 \hfill$\Box$

\subsubsection{Proof of Theorem 2.3.}
Consider the case with the nonparametric estimator $\hat p(X_i)$. 
It is well known that, by Assumptions~7 and 8,
$$\mathop{\sup}_{x\in\chi}\mid \hat
p(X) -p(X)  \mid=O_p(h_0^{{s}}+\sqrt{\frac {\log n}{ {nh_0^k}}})=o_p(h_0^{{s/2}}), \, s\geq k+2\geq 4.$$
It follows that, from the proof in (\ref{theo2}), and under Assumptions~7, 8 and 9 on the propensity score and the bandwidths,
\begin{eqnarray}\label{cn21}
\begin{split}
|C_{n,2}|&=\Big |\sum_{i=1}^n
\left(\frac{D_i}{ {\hat p}(X_i)}
-\frac{D_i}{ {
p}(X_i)}\right)K_iR_{n,i}(Y_i,\theta_\tau) \Big |\\
&\leq C\mathop{\sup}_{x} \mid \hat p(X) -p(X)  \mid
\sum_{i=1}^n  \mid K_iR_{n,i}(Y_i,\theta_\tau) \mid \\
&=O_p(n(h_0^{{s/2}})
O_p(n^{-3/4}h^{-1/4})=O_p((n/h)^{1/4}(h_0^{{s/2}})=o_p(1).
\end{split}
\end{eqnarray}

Next we want to get the asymptotic expansion for $C_{n,1}$
taking account of the nonparametric estimator $\hat p(X).$

Define the first and second partial derivatives of
$\phi(p(X_i),Z_i)$ w.r.t. the argument $p$ as
$\phi_p(p(X_i),Z_i)$ and $\phi_{pp}(p(X_i),Z_i)$,
respectively. By  Taylor expansion around $p(X_i)$, we have
\begin{align}\label{cn1}
C_{n,1}&=\frac{1}{\sqrt{nh}}
\sum_{i=1}^n\phi_p(p(X_i),Z_i)K_i(\hat p(X_i)-p(X_i))+
\frac{1}{\sqrt{nh}}
\sum_{i=1}^n\phi_{pp}(\tilde p(X_i),Z_i)K_i(\hat
p(X_i)-p(X_i))^2 \nonumber\\
&=\frac{1}{\sqrt{nh}}
\sum_{i=1}^nS_p(X_i)K_i(\hat
p(X_i)-p(X_i))+\frac{1}{\sqrt{nh}}
\sum_{i=1}^n\xi_iK_i(\hat p(X_i)-p(X_i)) \nonumber\\
&\quad+
\frac{1}{2\sqrt{nh}}
\sum_{i=1}^n\phi_{pp}( \tilde p(X_i),Z_i)K_i(\hat
p(X_i)-p(X_i))^2\nonumber \\
&:=J_1+J_2+J_3,
\end{align}
where  $S_p(X_i)=E(\phi_p(p(X_i),Z_i)\mid X_i)$,
$\xi_i=\phi_p(p(X_i),Z_i)-E(\phi_p(p(X_i),Z_i)\mid X_i)$ is
$i.i.d.$ and have conditional mean zero given $x_i$ and $\tilde
p(X_i)$ is a value between $p(X_i)$ and $\hat p(X_i)$.

Recall that $\hat p(X_i)=\sum_{j\neq i}w_{ij}D_j,$ where
$w_{ij}=\frac{L\left(\frac{X_j-X_i}{h_0}\right)}{\sum_{j\neq
i}L\left(\frac{X_j-X_i}{h_0}\right)}.$ We can decompose $\hat
p(X_i)$ as
\begin{eqnarray}
\hat p(X_i)-p(X_i)=\sum_{j\neq
i}w_{ij}\epsilon_j+\vartheta_n(X_i),
\end{eqnarray}
where $\epsilon_j=D_j-p(X_j),$ $\vartheta_n(X_i)=\sum_{j\neq
i}w_{ij}p(X_j)-p(X_i).$

Based on this formula, we can further write
\begin{eqnarray}\label{j1}
&&\quad\quad \quad J_1=\frac{1}{\sqrt{nh}}\sum_{i=1}^{n}S_p(X_i) K_i(\hat
p(X_i)- p(X_i))\nonumber\\
&&\quad\quad \quad =\frac{1}{\sqrt{nh}}\sum_{j=1}^{n}S_p(X_j)
K_j\epsilon_j+\frac{1}{\sqrt{nh}}\sum_{j=1}^n\epsilon_j
\left[\sum_{i\neq j}w_{ij}S_p(X_i)K_i-S_p(X_j)K_j
\right]\nonumber\\
&&\qquad\qquad+\frac{1}{\sqrt{nh}}\sum_{i=1}^{n}S_p(X_i)
K_i\vartheta_n(X_i)\nonumber\\
&&\quad\quad\qquad:=J_{11}+J_{12}+J_{13}.
\end{eqnarray}
Here $J_{11}=\frac{1}{\sqrt{nh}}\sum_{j=1}^{n}S_p(X_j)
K_j\epsilon_j=O_p(1)$ is the sum of $i.i.d.$ variables and is
serving as the leading term of $C_{n,1}.$ In the following, we aim to show  the
rest terms $J_{12}$, $J_{13}$, $J_2$  and $J_{3}$ are all
$o_p(1)$ such that $C_{n,1}=J_{11}+o_p(1)$.

{ First, we deal with $J_{12}$}. Clearly, its mean is zero.
Further, for each $j$, $\left[\sum_{i\neq
j}w_{ij}S_p(X_i)K_i-S_p(X_j)K_j \right] $ is rather like
$\left[\sum_{i\neq j}w_{ji}S_p(X_i)K_i-S_p(X_j)K_j \right] $
in terms of its magnitude. See \cite{ichimura2005} and
\cite{linton2001}. We now focus on
$\left[\sum_{i\neq j}w_{ji}S_p(X_i)K_i-S_p(X_j)K_j \right]$,
which can be regarded as the bias function from smoothing
$S_p(X_j)K_j$.

Since $X=(X_1,X_{(2)}),
L(\frac{X-X_j}{h_0})=
L_1(\frac{X_1-X_{1j}}{h_0})*
L_2(\frac{X_2-X_{2j}}{h_0}),$
we have, noting that $\hat f-f=o_p(1)$, and the kernel function $K(\cdot)$ is $s^*\ge s$ times continuously differentiable,
 \begin{eqnarray*}
& &E(\sum_{j\not =i, i=1}^nw_{ji}S_p(X_i)K_i\mid X_j)\\
&&=\frac{1+o_p(1)}{h_0^kf(X_j)}\int
L_1(\frac{u_{1i}-X_{1j}}{h_0})*L_2
(\frac{u_{2i}-X_{2j}}{h_0})K(\frac{u_{1i}
-x_{10}}{h})S_p(u_i)f(u_i)du\\
& &=\frac{1+o_p(1)}{f(X_j)}\int L_1(v_1)L_2(v_2)K(\frac{X_{1j}-
x_{10}}{h}+v_1\frac{h_0}{h})S_p(X_j+h_0v)f(X_j+h_0v)dv\\
& &=S_p(X_j)K_j+O_p(\frac{h_0^s}{h^s}).
\end{eqnarray*}
Thus, $\sup_j\mid \sum_{i \neq
j}w_{ji}S_p(X_i)K_i-S_p(X_j)K_j \mid=O_p(\frac{h_0^s}{h^s}).$
By Assumption 9 with $h_0^{2s}h^{-(2s+1)}\rightarrow0$, we have
\begin{eqnarray*}
\sup_i\mid \frac{1}{\sqrt{h}}(\sum_{i \neq
j}w_{ij}(S_p(X_i)K_i-S_p(X_j)K_j) \mid=O_p(\frac {h_0^{2s}}{h^{(2s+1)}})=o_p(1).
\end{eqnarray*}
Since $\epsilon_j=D_j-E(D_j\mid X_j)$, $\epsilon_j$' s are
mutually independent conditional on the sample path of the
$X_j$' s, thus we have $J_{12}=o_p(1).$

{ {Next, we  prove $J_{13}=o_p(1).$ }} Based on $Lemma$
$A.1(b)$ in \cite{abrevaya2015}, noting that it is actually the bias term of a nonparametric estimation,
$$\sup_{x\in \chi}\mid
\vartheta_n(x) \mid=O_p(h_0^s). $$
 Thus we have, by Assumptions 8 and 9 with $nhh_0^{2s}\rightarrow 0$,
\begin{eqnarray*}
\mid J_{13}\mid &=&\mid\frac{1}{\sqrt{nh}}\sum_{i=1}^{n}
S_p(X_i)K_i\vartheta_n(x_i)\mid \nonumber\\
&\leq & \sqrt{nh} sup_i\mid \vartheta_n(x_i)\mid \frac{1}{nh}
\sum_{i=1}^{n}\mid S_p(X_i)\mid \mid K_i\mid \nonumber\\
&=&\sqrt{nh} O_p(h_0^s).O_p(1)=o_p(1).
\end{eqnarray*}

As for $J_2$, we have, by Assumptions 8 and 9 again,
\begin{eqnarray*}
\frac{1}{\sqrt{h}}\mid K_i (\hat p(X_i)-p(X_i))\mid\leq
C\frac{1}{\sqrt{h}}\sup_i\mid\hat p(X_i)-p(X_i)\mid
= O_p(h^{-1/2}h_0^{{s/2}}).
\end{eqnarray*}
since $\xi_i$'s are mutually independent conditional on the
sample path of $x_i$'s, the similar argument for proving
$J_{12}$ can be applied to derive $J_2=O_p(h^{-1/2}h_0^{{s/2}})=o_p(1).$

When it comes to $J_3$, since $\tilde p(X_i)$ is between $\hat
p(X_i)$ and $p(X_i)$, $0<\mid \tilde p(X_i)\mid <1.$ Furthermore,
similarly as the previous discussion, we can show that $sup_x\mid
\hat p(X)-p(X)\mid^2=O_p((h_0^{{s}}+\sqrt{\frac {\log n}{ {nh_0^k}}})^2)$ and
$\frac{1}{nh}\sum_{i=1}^{n}\mid
\psi_{pp}(\tilde p(X_i),\omega_i)K_i \mid=O_p(1). $ Thus, Assumptions~8 and 9 implies that $\mid
J_{3}\mid\leq
\sqrt{nh}O_p((h_0^{{s}}+\sqrt{\frac {\log n}{ {nh_0^k}}})^2)O_p(1)=O_p(\sqrt{nh}h_0^{{s}})=o_p(1).$

Until now, we have proved under the nonparametric estimator
$\hat p(X)$,
\begin{center}
$C_{n,2}=o_p(1);$
$C_{n,1}=J_{11}+o_p(1)=\frac{1}{\sqrt{nh}}
\sum_{j=1}^{n}S_p(X_j) K_j\epsilon_j+o_p(1).$
\end{center}
Thus we have
\begin{eqnarray}
\begin{split}
A_{n,1}&=\theta_\tau^{\top}J_{11}+o_p(1)=
\theta_\tau^{\top}\frac{1}{\sqrt{nh}}\sum_{i=1}^{n}S_p(X_i)
K_i\epsilon_i+o_p(1) \nonumber\\
&=u_\tau\frac{1}{\sqrt{nh}}\sum_{i=1}^{n}s^1_p(X_i)\epsilon_i
K_i+v_\tau^{\top}\frac{1}
{\sqrt{nh}}\sum_{i=1}^{n}S^{(2)}_p(X_i)\epsilon_i
K_i+o_p(1).
\end{split}
\end{eqnarray}

Note that from the proof of Lemma~2.1, we have
\begin{eqnarray*}
&A_{n,2}=u_\tau \frac{1}{\sqrt{nh}}\sum_{i=1}^{n}\frac{D_i}{p(X_i)}
(\mathbb I(Y_i\leq
\tilde q_{1,\tau}(X_{1i}))-\tau)K_i+\frac{1}{2}u_\tau^2f_1( q_{1,\tau}(x_{10})\mid x_{10})f(x_{10})\\
&+v_\tau^{\top}W_{n2}+
\frac{1}{2}f_1( q_{1,\tau}(x_{10})\mid x_{10})f(x_{10})v_\tau^{\top}\mu_2(K)v_\tau+o_p(1).
\end{eqnarray*}
That implies
\begin{eqnarray*}
\begin{split}
&A_{n,1}+A_{n,2}\\
&=u_\tau \frac{1}{\sqrt{nh}}\sum_{i=1}^{n}\frac{D_i}{p(X_i)}
(\mathbb I(Y_i\leq
\tilde q_{1,\tau}(X_{1i}))-\tau)K_i+\frac{1}{2}u_\tau^2f_1( q_{1,\tau}(x_{10})\mid x_{10})f(x_{10}) \\
&\quad +v_\tau^{\top}\frac{1}
{\sqrt{nh}}\sum_{i=1}^{n}S^{(2)}_p(X_i)\epsilon_iK_i+
\frac{1}{2}f_1( q_{1,\tau}(x_{10})\mid x_{10})f(x_{10})v_\tau^{\top}\mu_2(K)v_\tau+o_p(1),
\end{split}
\end{eqnarray*}
where $S_p(X_i)=E(\phi_p(p(X_i),Z_i)\mid X_i)=\left(
\begin{array}{c}
 -\frac{E(\mathbb I(Y_i(1)^*\leq0)-\tau\mid X_i)}{p(X_i)} \\
-\frac{E(\mathbb I(Y_i(1)^*\leq0)-\tau\mid X_i)}{p(X_i)}\frac{X_{1i}-x_{10}}{h} \\
 \end{array}
\right):=\left(
           \begin{array}{c}
             s_p^1(X_i) \\
             S_p^{(2)}(X_i)\\
           \end{array}
         \right)$.

It follows that
\begin{eqnarray*}
\begin{split}
&A_{n,1}+A_{n,2}\\
&=u_\tau \frac{1}{\sqrt{nh}}\sum_{i=1}^{n}
\frac{D_i(I(Y_i\leq \tilde q_{1,\tau}(X_{1i}))-\tau)-E(\mathbb I(Y_i(1)\leq
\tilde q_{1,\tau}(X_{1i}))-\tau\mid X_i)\epsilon_i}{p(X_i)}
K_i\\
&\quad +\frac{1}{2}u_\tau^2f_1( q_{1,\tau}(x_{10})\mid x_{10})f(x_{10})
+v_\tau^{\top}\frac{1}
{\sqrt{nh}}\sum_{i=1}^{n}S^{(2)}_p(X_i)\epsilon_iK_i\\
&\quad +
\frac{1}{2}f_1( q_{1,\tau}(x_{10})\mid x_{10})f(x_{10})v_\tau^{\top}\mu_2(K)v_\tau+o_p(1).
\end{split}
\end{eqnarray*}
Since $A_{n,3}$ does not depend on $\theta_\tau$ and based on
the convexity lemma \citep{pollard1991}, we can get the
asymptotically linear representation of $\hat\Delta^{ncqte}_\tau(x_{10})$
under the nonparametric estimator $\hat p(X)$ as
\begin{eqnarray*}
\sqrt{nh}\left(\hat q^{ncqte}_{1,\tau }(x_{10})- {q_{1,\tau}
(x_{10})}\right)=-\frac{1}{\sqrt{nh}}
\frac{1}{f(x_{10})}\sum_{i=1}^{n}
\left(\frac{D_i}{p(X_i)} \eta_{1,
\tau}(Y_i)+n_p^1(X_i)\epsilon_i\right)K_i+o_p(1),
\end{eqnarray*}
where $Z_i=(X_i,D_i,Y_i),K_i=K(\frac{X_{1i}-x_{10}}{h})$,
$\eta_{1,\tau}(Y_i)=\frac{\mathbb I(Y_i\leq
\tilde q_{1,\tau}(X_{1i}))-\tau}
{f_1(q_{1,\tau}(x_{10})\mid x_{10})}$,
$n_p^1(X_i)=\frac{-E(\mathbb I(Y_i(1)\leq
\tilde q_{1,\tau}(X_{1i}))-\tau\mid X_i)}{p(X_i)f_1(q_{1,\tau}(x_{10})\mid x_{10})}$.

By the analogous proof for $\hat q^{ncqte}_{1,\tau}(x_{10})$, we can get
that the asymptotically linear representation of $\hat
q_{0,\tau}(x_{10})$ as follows
\begin{eqnarray*}
& &\sqrt{nh}\left(\hat q^{ncqte}_{0,\tau }(x_{10})- {q_{0,\tau}
(x_{10})}\right)=-\frac{1}{\sqrt{nh}}
\frac{1}{f(x_{10})}\sum_{i=1}^{n}\left(
\frac{1-D_i}{1-p(X_i)} \eta_{0,
\tau}(Y_i)+n_p^0(X_i)\epsilon_i\right)K_i+o_p(1),
\end{eqnarray*}
where $\eta_{0,\tau}(Y_i)=\frac{\mathbb I(Y_i\leq
\tilde q_{0,\tau}(X_{1i}))-\tau}{f_0(q_{0,\tau}
(x_{10})\mid x_{10})}$ and  $n_p^0(x_i)=\frac{E(\mathbb I(Y_i(0)\leq
\tilde q_{0,\tau}(X_{1i}))-\tau\mid X_i)}{(1-p(X_i))f_0(q_{0,\tau}(x_{10})\mid x_{10})}$.

Thus we can get the asymptotically linear representation of
$\hat\Delta^{ncqte}_\tau(x_{10})$ as
\begin{eqnarray*}
& &\sqrt{nh}\left(\hat \Delta^{ncqte}_{\tau }(x_{10})-
{\Delta_{\tau} (x_{10})}\right)=-\frac{1}{\sqrt{nh}}
\frac{1}{f(x_{10})}\sum_{i=1}^{n}\phi_2(p(X_i),Z_i)K_i+o_p(1),\\
\end{eqnarray*}
where $ \phi_2(p(X_i),Z_i)=\phi_1
(p(X_i),Z_i)-n_p(X_i)\epsilon_i$. Further, $\phi_1
(p(X_i),Z_i)=
\frac{D_i}{p(X_i)} \eta_{1,
\tau}(Y_i) -\frac{1-D_i}{1-p(X_i)} \eta_{0,\tau}(Y_i)$,
$ n_p(X_i)=-n_p^1(X_i)+n_p^0(X_i)
=\frac{m_{1,\tau}(X_i)}
{p(X_i)}+\frac{m_{0,\tau}
(X_i)}{(1-p(X_i))} .$

Again, we  analyse the asymptotic Bias and variance
of $\hat \Delta_\tau(x_{10}).$

{\it{Bias of $NCQTE(X_{10})$}.}
Since $E(\epsilon\mid X)=E(D-p(X)\mid X)=0$ and $K$ is a kernel of order $s_1$, we can get the the bias of $NCQTE(x_{10})$ similarly as
$OCQTE(x_{10})$, i.e.
\begin{eqnarray}
Bias(\hat
\Delta^{ncqte}_\tau(x_{10}))=O_p(\mu_{s_1}(K)h^{s_1}),
\end{eqnarray}
where $\mu_{s_1}(K)=\int u^{s_1}K(u)du.$

{\it{Variance of $NCQTE(X_{10})$}.}
Let $\sigma_{ncqte}^{*2}(X_1)=E(\phi_2^2(p(X),z)\mid X_1).$ Similarly
as the proof of $PCQTE(x_{10})$, we have
\begin{eqnarray*}
\begin{split}
Var(\hat\Delta^{ncqte}_{\tau}(x_{10}))=\frac{1}{nh}\frac{\parallel
K\parallel_2^2\sigma_{ncqte}^{*2}(x_{10})}{f(x_{10})}+o\left(\frac{1}{nh}\right),
\end{split}
\end{eqnarray*}
where $\parallel K\parallel_2^2=\int K^2(t)dt.$

Suppose that the assumptions 1 through 5 and 7
through 9 are satisfied for some $s_1\geq s\geq k+l$, we can get
 \begin{eqnarray*}
\sqrt{nh}\left(\hat {\Delta}^{ncqte}_\tau (x_{10})-\Delta_\tau
(x_{10})-O_p(\mu_{s_1}h^{s_1})\right)
\stackrel{\mathcal{D}}{\longrightarrow}N
\left(0,\frac{\parallel K \parallel_2^2\sigma_{ncqte}^{*2}(x_{10})}
{f(x_{10})}\right).
\end{eqnarray*}
The proof for $Theorem$~2.3 is concluded.  \hfill$\Box$

\subsubsection{Proof of Corollary~2.1.}
We show that $NCQTE(x_{10})$ can be more efficient
than $PCQTE(x_{10})$ i.e.
$\sigma_{ncqte}^{*2}(x_{10})\leq \sigma_{pcqte}^2(x_{10}).$

Note that $\phi_2(p(X),Z)=\phi_1
(p(X),Z)-n_p(X)\epsilon$. Thus
$\sigma_{ncqte}^{*2}(x_{10})=\sigma_{pcqte}^2(x_{10})+
E(n_p(X)^2\epsilon^2\mid
X_1=x_{10})-2E(\phi_1(p(X),Z)n_p(X)\epsilon\mid X_1=x_{10}).$ Further,  some elementary calculations yield that
 \begin{eqnarray}\label{ncqte_comvar}
E(n_p(X)^2\epsilon^2\mid
X_1=x_{10})=E(\phi_1(p(X),Z)n_p(X)\epsilon\mid X_1=x_{10})\nonumber\\
=E\left[p(X)(1-p(X))\left(\frac{m_{1,\tau}(X)}
{p(X)}+\frac{m_{0,\tau}(X)}{1-p(X)}\right)^2\mid
X_1=x_{10}\right].
 \end{eqnarray}
 Thus we have
 $\sigma_{pcqte}^{2}(x_{10})=\sigma_{ncqte}^{*2}(x_{10})+E(n_p(X)^2\epsilon^2\mid
 x_{10})\geq\sigma_{ncqte}^{*2}(x_{10}).$
The proof of  $Corollary$~2.1 is completed.
 \hfill$\Box$

\subsubsection{Proof of Theorem~2.4.}
\renewcommand{\theequation}{S.\arabic{equation}}
\setcounter{equation}{0}
\renewcommand{\thelemma}{S.\arabic{lemma}}
\setcounter{lemma}{0}

\quad  Recall that the seimiparametric estimator of
 $p(X_i)=p(\alpha^{\top}X_i)$ is
\begin{eqnarray*}
\hat p(\hat\alpha^{\top} X_i)=\frac{\frac{1}{nh_2^q}\sum_{j\neq
i}D_jH\left(\frac{\hat\alpha ^{\top}X_j-\hat\alpha
^{\top}X_i}{h_2}\right)}{\frac{1}{nh_2^q}\sum_{j\neq
i}H\left(\frac{\hat\alpha ^{\top}X_j-\hat\alpha
^{\top}X_i}{h_2}\right)}\label{semi+equ2},
\end{eqnarray*}
where $\hat \alpha -\alpha=O_p(n^{-1/2}).$
The difference between $\hat p(\hat \alpha^\top X_i)$ and
$\hat p(\alpha^\top X_i)$ comes from the difference between
$H\left(\frac{\hat\alpha ^{\top}X_j-\hat\alpha
^{\top}X_i}{h_2}\right)$ and $H\left(\frac{\alpha
^{\top}X_j-\alpha ^{\top}X_i}{h_2}\right)$. Further, we can
rewrite the former as
\begin{eqnarray*}
H\left(\frac{\hat\alpha ^{\top}X_j-\hat\alpha
^{\top}X_i}{h_2}\right)=H\left(\frac{\alpha ^{\top}X_j-\alpha
^{\top}X_i}{h_2}+\frac{(\hat\alpha-\alpha)
^{\top}(X_j-X_i)}{h_2}\right).
\end{eqnarray*}

Thus we have
$$M(\hat \alpha,X_i )=M(\alpha,X_i)+M_{1,n}(\alpha,X_i)(\hat
\alpha-\alpha)+o_p(M_{1,n}(\alpha,X_i)(\hat\alpha-\alpha)),$$
where $M(\hat \alpha,X_i)=\frac{1}{nh_2^q}\sum_{j\neq
i}H\left(\frac{\hat\alpha ^{\top}X_j-\hat\alpha
^{\top}X_i}{h_2}\right)$, $M_{1,n}(\alpha,X_i)=\frac{1}{nh_2^{q+1}}\sum_{j\neq
i}H'\left(\frac{\alpha ^{\top}X_j-\alpha
^{\top}X_i}{h_2}\right)\\(X_j-X_i)^{\top}.$

Using the facts $\int H'(u)du=0$ and $\int uH'(u)du=-1$, we have
$M_{1,n}(\alpha,X_i)=O_p(1).$ It follows that
$M(\hat\alpha,X_i)=M(\alpha,X_i)+O_p(\frac{1}{\sqrt{n}}).$
Further, it is not hard to show that $\hat
p(\hat\alpha^{\top}X_i)=\hat p(\alpha^{\top}X_i)
+O_p(\frac{1}{\sqrt{n}})$ uniformly for all $i$.

Based on this result, It follows that
\begin{eqnarray*}
\begin{split}
&\frac{1}{n}\sum_{i=1}^{n}\frac{D_i}{
{\hat p}(\hat \alpha^{\top}x_i)}\rho_{\tau}(Y_i-a-b(X_{1i}-x_{10}))K_i\\
&=\frac{1}{n}\sum_{i=1}^{n}\frac{D_i}{
{\hat p}(\alpha^{\top}x_i)}\rho_{\tau}(Y_i-a-b(X_{1i}-x_{10}))K_i+O_p(\frac{1}{\sqrt{n}})\\
&=\frac{1}{n}\sum_{i=1}^{n}\frac{D_i}{
{\hat p}(\alpha^{\top}x_i)}\rho_{\tau}(Y_i-a-b(X_{1i}-x_{10}))K_i+o_p(\frac{1}{\sqrt{nh}})
\end{split}
\end{eqnarray*}

That is, the $SCQTE(x_{10})$, based on $\hat p(\hat \alpha^{\top}X)$, is asymptotically as efficient as that based on $\hat p(\alpha^{\top}X)$. In other words, for the asymptotic efficiency, the estimator $\hat \alpha $ plays no role. Thus we can simply, without confusion, regard $\hat \alpha$ as $\alpha$ and rewrite $G_\tau(a,b,\hat p)$ as
\begin{eqnarray*}
&&G_{\tau}( a,b,\hat p)\\
&&=\sum_{i=1}^{n}\left(\frac{D_i}{ {\hat
p}(\alpha^{\top}X_i)}-\frac{D_i}{ {
p}(\alpha^{\top}X_i)}\right)
\left[\rho_{\tau}(Y_i-a-b(X_{1i}-x_{10}))-
\rho_{\tau}(Y_i-\tilde q_{1,\tau}(X_{1i}))\right]K_i\nonumber\\
\end{eqnarray*}

\begin{eqnarray*}
&&+\sum_{i=1}^n\frac{D_i}{{
p}(\alpha^{\top}X_i)}\left[\rho_{\tau}(Y_i-a-b(X_{1i}-x_{10}))
-\rho_{\tau}(Y_i-\tilde q_{1,\tau}(X_{1i}))\right]K_i\nonumber\\
&&+\sum_{i=1}^n\frac{D_i}{{ \hat
p}(\alpha^{\top}X_i)}\rho_{\tau}(Y_i-\tilde q_{1,\tau}(X_{1i}))K_i\nonumber\\
&&:=A_{n,1}+A_{n,2}+A_{n,3}.
\end{eqnarray*}
From the proof of Lemma~2.1, it can also be verified that
\begin{eqnarray*}
&A_{n,2}=u_\tau \frac{1}{\sqrt{nh}}\sum_{i=1}^{n}\frac{D_i}{p(\alpha^{\top}x_i)}
(\mathbb I(Y_i\leq
\tilde q_{1,\tau}(X_{1i}))-\tau)K_i+\frac{1}{2}u_\tau^2f_1( q_{1,\tau}(x_{10})\mid x_{10})f(x_{10})\\
&+v_\tau^{\top}W_{n2}+
\frac{1}{2}f_1( q_{1,\tau}(x_{10})\mid x_{10})f(x_{10})v_\tau^{\top}\mu_2(K)v_\tau+o_p(1).
\end{eqnarray*}

Next we study the asymptotic behavior of $A_{n,1}$. Again, similarly rewrite $A_{n,1}$ as
\begin{eqnarray}
\begin{split}
&A_{n,1}=\sum_{i=1}^{n}\left(\frac{D_i}{ {\hat
p}(\alpha^{\top}X_i)}-\frac{D_i}{ {
p}(\alpha^{\top}X_i)}\right)[\frac{\theta_\tau^{\top}X_{1i}^*}{\sqrt{nh}}\eta_i^*+R_{n,i}(Y_i,\theta_\tau)]
K_i\nonumber\\
&=\frac{\theta_\tau^{\top}}{\sqrt{nh}}
\sum_{i=1}^n[\phi(Z_i,\hat
p(\alpha^{\top}X_i))-\phi(Z_i,p(\alpha^{\top}X_i))]K_i+\tilde C_{n,2}:=\theta_\tau^{\top}
\tilde C_{n,1}+\tilde C_{n,2}.
\end{split}
\end{eqnarray}
Here $\phi(p(\alpha^{\top}X_i),Z_i)
=\frac{D_iX_{1i}^*}{p(\alpha^{\top}X_i)}
(\mathbb I\left(Y_i^*\leq
0\right)-\tau),\tilde C_{n,1}=\frac{1}{\sqrt{nh}}
\sum_{i=1}^n[\phi(Z_i,\hat p(\alpha^{\top}X_i))-\phi(Z_i,p(\alpha^{\top}X_i))]K_i,$
and
$\tilde C_{n,2}=\sum_{i=1}^n\left(\frac{D_i}{ {\hat
p}(\alpha^{\top}X_i)}-\frac{D_i}{ {
p}(\alpha^{\top}X_i)}\right)R_{n,i}(Y_i,\theta_\tau)K_i.$

Similarly as $C_{n,2}$ in the case of $NCQTE(x_{10})$, we can show that $\tilde C_{n,2}=o_p((n/h)^{1/4}(h_2^{{s_2/2}})=o_p(1)$ by Assumptions~7', 8', and 9'.

 By a Taylor expansion around $p(\alpha^{\top}X_i)$, we can also decompose $\tilde C_{n,1}$ as
\begin{eqnarray*}
\begin{split}
\tilde C_{n,1}&=\frac{1}{\sqrt{nh}}
\sum_{i=1}^n\phi_p(p(\alpha^{\top}X_i),Z_i)K_i(\hat p(\alpha^{\top}X_i)-p(\alpha^{\top}X_i))\\
& \quad +
\frac{1}{\sqrt{nh}}
\sum_{i=1}^n\phi_{pp}(\tilde p(\alpha^{\top}X_i),Z_i)K_i(\hat
p(\alpha^{\top}X_i)-p(\alpha^{\top}X_i))^2\\
&=\frac{1}{\sqrt{nh}}
\sum_{i=1}^ns_p(X_{1i},\alpha^{\top}X_i)K_i(\hat
p(\alpha^{\top}X_i)-p(\alpha^{\top}X_i))+\frac{1}{\sqrt{nh}}
\sum_{i=1}^n\xi_iK_i(\hat p(\alpha^{\top}X_i)-p(\alpha^{\top}X_i))\\
&\quad+
\frac{1}{2\sqrt{nh}}
\sum_{i=1}^n\phi_{pp}(\tilde p(\alpha^{\top}X_i),Z_i)K_i(\hat
p(\alpha^{\top}X_i)-p(\alpha^{\top}X_i))^2\\
&:=\tilde J_1+\tilde J_2+\tilde J_3,
\end{split}
\end{eqnarray*}
where $s_p(X_{1j},\alpha^{\top}X_j)=E(\phi_p(p(\alpha^{\top}X_i),Z_i)\mid X_{1j},\alpha^{\top}X_j),\xi_i=\phi_p(p(\alpha^{\top}X_i),Z_i)-E(\phi_p(p(\alpha^{\top}X_i),Z_i)\mid X_{1j},\alpha^{\top}X_j).$

Further, note that $\hat p(\alpha^{\top}X_i)=\sum_{j:j\neq i}^nw_{ij}\varepsilon_j+\beta_n(\alpha^{\top}X_i)$, $\varepsilon_j:=D_j-p(\alpha^{\top}X_j), w_{ij}=\frac{H((\alpha^{\top}X_j-\alpha^{\top}X_i)/h_2)D_j}{\sum_{j:j\neq i}H((\alpha^{\top}X_j-\alpha^{\top}X_i)/h_2)}.$ Thus we have
\begin{eqnarray}\label{tj1}
\tilde J_1&=&\frac{1}{\sqrt{nh}}\sum_{i=1}^ns_p(X_{1i},\alpha^{\top}X_i)(\hat p(\alpha^{\top}X_i)-p(\alpha^{\top}X_i))K(\frac{X_{1i}-x_{10}}{h})\nonumber\\
&=&\frac{1}{\sqrt{nh}}\sum_{i=1}^ns_p(X_{1i},\alpha^{\top}X_i)K(\frac{X_{1i}-x_{10}}{h})\sum_{j\neq i}^nw_{ij}\varepsilon_j\nonumber\\
&&+\frac{1}{\sqrt{nh}}\sum_{i=1}^ns_p(X_{1i},\alpha^{\top}X_i)
\beta_n(\alpha^{\top}X_i)K(\frac{X_{1i}-x_{10}}{h})\nonumber\\
&:=&\tilde J_{11}+\tilde J_{13}.
\end{eqnarray}
{
 As the  proofs for  $\tilde J_{13}=o_p(1)$, $\tilde J_{2}=o_p(1)$ and $\tilde J_{3}=o_p(1)$ under Assumptions~ 7', 8' and 9' are very similar to those in the nonparametric setting in Theorem~2.3, we then omit the details here. However, as for $\tilde J_{11}$, its asymptotic behavior varies with whether $X_1$ is a subset of $\alpha^{\top}X$ or not, which
can be summarized as the following lemmas. The proofs of these lemmas will be
presented in next subsection. }

\begin{lemma}{\label{out_lemma}}
If $X_1\sqsubset^{l-r} \alpha^{\top}X$ with $s_2(2-l/r)+l>0$ and $0< r\leq l$, under the assumptions 1 through 4 and 7'
through 9' with $s^*\geq s_2\geq 2$, we can obtain \begin{eqnarray} \tilde J_{11}=o_p(1).
\end{eqnarray}
Further, the asymptotically linear representation of $\hat
\Delta^{scqte}_{\tau}(x_{10})$ can be
 \begin{eqnarray*}
& &\sqrt{nh^l}\left(\hat \Delta^{scqte}_{\tau }(x_{10})-
{\Delta_{\tau} (x_{10})}\right)=-\frac{1}{\sqrt{nh^l}}
\frac{1}{f(x_{10})}\sum_{i=1}^{n}\phi_1(p(\alpha^\top
X_i),Z_i)K_i+o_p(1).
\end{eqnarray*}
\end{lemma}

 {
\begin{lemma}{\label{in_lemma}}
If $X_1\subsetneq\alpha^{\top}X$, under the assumptions 1 through 4 and 7'
through 9' with $s^*\geq s_2\geq 2$, we can obtain
\begin{eqnarray}
\tilde J_{11}=\frac{1}{\sqrt{nh}}\sum_{i=1}^ns_p(\alpha^{\top}X_i)
\epsilon_iK_{1i}+o_p(1).
\end{eqnarray}
The corresponding asymptotically linear representation of $\hat \Delta^{scqte}_\tau(x_{10})$ is
 \begin{eqnarray*}
& &\sqrt{nh^l}\left(\hat \Delta^{scqte}_{\tau }(x_{10})-
{\Delta_{\tau} (x_{10})}\right)=-\frac{1}{\sqrt{nh^l}}
\frac{1}{f(x_{10})}\sum_{i=1}^{n}\phi_3^*(p(\alpha^\top
X_i),Z_i)K_i+o_p(1),
\end{eqnarray*}
where $\phi^*_3(p(\alpha^{\top}X_i),Z_i)=\psi(p
(\alpha^{\top}X_i,Z_i))-e_p(\alpha^{\top}X_i)\epsilon_i,$
$e_p(\alpha^\top x_i)=\frac{m_{1,\tau}(\alpha^\top X_i)}
{p(\alpha^\top X)}+\frac{m_{0,\tau}
(\alpha^\top X_i)}{1-p(\alpha^\top X)}$, $m_{j,\tau}(\alpha^\top X_i)=\frac{E[\mathbb
I(Y_i(j)\leq \tilde q_{j,\tau}(X_{1i}))-\tau\mid
\alpha^{\top}X_i]}{f_j(q_{j,\tau}(x_10)\mid x_{10})},j=0,1,$
and $\epsilon=D-p(\alpha^\top X)$.
\end{lemma}
\bigskip
}

%
%

Based on these two lemmas, we can easily get the first two parts of Theorem~2.4 by the analogue aforementioned discussion. Finally, we aim to show that $\sigma^2_{scqte}(x_{10})\geq  \sigma_{scqte}^{2*}(x_{10}),$
where  $
\sigma_{scqte}^2(x_{10})=E(\phi_1^2(p(X),Z)\mid X_1=x_{10}), ~ \sigma_{scqte}^{*2}(x_{10})=E(\phi_3^{*2}(p(X),Z)\mid X_1=x_{10}).$

 Note that $\phi_3^*(p(\alpha^{\top}X_i),Z_i)=\phi_1
(p(\alpha^{\top}X_i),Z_i)-e_p(\alpha^{\top}X_i)\epsilon_i$. Thus, following the same arguments of (\ref{ncqte_comvar}),  \begin{eqnarray}\label{ncqte_comvar}
E(e_p(\alpha^{\top}X)^2\epsilon^2\mid
X_1=x_{10})=E(\phi_1(p(\alpha^{\top}X),Z)e_p(\alpha^{\top}X)\epsilon\mid X_1=x_{10})\nonumber\\
=E\left[p(\alpha^{\top}X)(1-p(\alpha^{\top}X))\left(\frac{m_{1,\tau}(\alpha^{\top}X)}
{p(\alpha^{\top}X)}+\frac{m_{0,\tau}(\alpha^{\top}X)}{1-p(\alpha^{\top}X)}\right)^2\mid
X_1=x_{10}\right].
 \end{eqnarray}
 Thus we can conclude the proof
of Theorem~2.4 immediately.
\begin{flushright}
  $\square$
\end{flushright}

 \subsubsection{Proof of Lemma~S.1}

Note that we assume $X_1\in R$, i.e. $l=1$ at the beginning, thus $X_1\sqsubset^{l-r}\alpha^{\top}X$ is equivalent to $X_1 \sqsubset^{0}\alpha^{\top}X$. Under this assumption, we can rewrite $\tilde J_{11}$ in the form of $U-$statistic. Denote
$\hat f(\alpha^{\top}x_i)=\sum_{j:j\neq i}H((\alpha^{\top}X_j-\alpha^{\top}X_i)/h_2),$ $ B_{ij}=s_p(X_{1i},\alpha^{\top}X_i)K(\frac{X_{1i}-x_{10}}{h})\frac{1}{f(\alpha^{\top}X_i)}\varepsilon_j,$ and
\begin{eqnarray*}
\tilde J^*_{11}=\frac{n-1}{\sqrt{nh}}\frac{1}{n(n-1)}\sum_{i=1}^n\sum_{j\neq i}^n[\frac{B_{ij}+B_{ji}}{2}]\frac{1}{h_2^q}H(\frac{\alpha^{\top}X_j-\alpha^{\top}X_i}{h_2}).
\end{eqnarray*}
We have
\begin{eqnarray*}
\mid \tilde J_{11}-\tilde J^*_{11}\mid \leq \sup_{\alpha^{\top}X\in\bigwedge^*}|\hat f(\alpha^{\top}X)-f(\alpha^{\top}X)| \inf_{\alpha^{\top}X\in\bigwedge^*} (\hat f(\alpha^{\top}X))^{-1}\mid \tilde J^*_{11}\mid.
\end{eqnarray*}
Here $\Lambda^*$ is a compact support of $\alpha^{\top}X$.
Since $h_2\rightarrow 0, nh_2^q\rightarrow \infty$, Assumptions~ 7', 8 and 9' imply that $\sup_{\alpha^{\top}X\in\bigwedge^*}|\hat f(\alpha^{\top}X)-f(\alpha^{\top}X)|=O_p(h_2^{s_2/2})$ and $ \hat f(\alpha^{\top}X)  $ is uniformly bounded, it follows that
\begin{eqnarray}
\mid \tilde J_{11}-\tilde J^*_{11}\mid \leq O_p(h_2^{s_2/2})\mid \tilde J^*_{11}\mid=o_p(\mid \tilde J^*_{11}\mid).
\end{eqnarray}

Note that $\tilde J^*_{11}$ is an $U-$stastistic.
Thus, we firstly compute the conditional expectation of $[\frac{B_{ij}+B_{ji}}{2}]\frac{1}{h_2^q}H(\frac{\alpha^{\top}X_j-\alpha^{\top}X_i}{h_2}).$
Since $E(\varepsilon_j\mid X_j)=0,$ it can be derived that \begin{center}
$E(\frac{1}{h_2^q}B_{ij}H(\frac{\alpha^{\top}X_j-
\alpha^{\top}X_i}{h_2})\mid X_i,D_i,Y_i)=0$ and $E(\frac{1}{h_2^q}B_{ij}H(\frac{\alpha^{\top}X_j-
\alpha^{\top}X_i}{h_2}))=0$.
\end{center}
Further we have
\begin{eqnarray*}
\begin{split}
&E\left(B_{ji}\frac{1}{h_2^q}H(\frac{\alpha^{\top}X_j-\alpha^{\top}X_i}{h_2})\mid X_i,D_i,Y_i\right)\\
&=E\left(s_p(X_{1j},\alpha^{\top}X_j)K(\frac{x_{1j}-x_{10}}{h})
\frac{1}{f(\alpha^{\top}X_j)}\frac{1}{h_2^q}H(\frac{\alpha^{\top}X_j-\alpha^{\top}X_i}{h_2})\mid X_i,D_i,Y_i\right)\varepsilon_i.\\
\end{split}
\end{eqnarray*}
Let $X_{1}=\nu_1,\alpha^{\top}X=\nu_2$, and denote $(\frac{\nu_1-\nu_{1i}}{h_2},\frac{\nu_2-\nu_{2i}}{h_2})$ as $(t_1,t_2)$. We have
\begin{eqnarray*}
\begin{split}
&E\left(s_p(X_{1j},\alpha^{\top}X_j)
K(\frac{X_{1j}-x_{10}}{h})
\frac{1}{f(\alpha^{\top}X_j)}\frac{1}{h_2^q}H(\frac{\alpha^{\top}X_j-\alpha^{\top}X_i}{h_2})\mid X_i,D_i,Y_i\right)\\
&=E\left(s_p(X_{1j},\alpha^{\top}X_j)
K(\frac{X_{1j}-x_{10}}{h})
\frac{1}{f(\alpha^{\top}X_j)}
\frac{1}{h_2^q}H(\frac{\alpha^{\top}X_j-\alpha^{\top}X_i}
{h_2})\mid X_i\right)\\
&=h_2\int K(\frac{t_1h_2}{h}+
\frac{\nu_{1i}-x_{10}}{h})\frac{1}
{f(\nu_{2i}+h_2t_2)}
H(t_2)s_p(\nu_{1i}+h_2t_1,\nu_{2i}+h_2t_2)\\
&\qquad\qquad\qquad\qquad\qquad \qquad\qquad\qquad\qquad\qquad \times f_{\nu_{12}}(\nu_{1i}+h_2t_1,\nu_{2i}+h_2t_2)
dt_1dt_2\\
&=h_2K(
\frac{\nu_{1i}-x_{10}}{h})\frac{1}
{f(\nu_{2i})}s_p(\nu_{1i},\nu_{2i})
f_{\nu_{12}}(\nu_{1i},\nu_{2i})\int H(t_2)
dt_1dt_2\\
&\qquad+\frac{h_2^{2}}{h}K'(
\frac{\nu_{1i}-x_{10}}{h})\frac{1}
{f(\nu_{2i})}s_p(\nu_{1i},\nu_{2i})
f_{\nu_{12}}(\nu_{1i},\nu_{2i})\int t_1H(t_2)
dt_1dt_2+o_p(\frac{h_2^{2}}{h}),
\end{split}
\end{eqnarray*}
where $f_{\nu_{12}}(\nu_1,\nu_2)$ is the joint density function of $(X_{1},\alpha^{\top}X).$

Under Assumptions $7'$ throught $9'$, we have
\begin{eqnarray*}
\begin{split}
&E\left(s_p(X_{1j},\alpha^{\top}X_j)K(\frac{X_{1j}-x_{10}}{h})
\frac{1}{f(\alpha^{\top}X_j)}\frac{1}{h_2^q}H(\frac{\alpha^{\top}X_j-\alpha^{\top}X_i}{h_2})\mid X_i,D_i,Y_i\right)\\
&=C_1h_2K(
\frac{X_{1i}-x_{10}}{h})\frac{1}
{f(\alpha^{\top}X_i)}s_p(X_{1i},\alpha^{\top}X_i))
f_{\nu_{12}}(X_{1i},\alpha^{\top}X_i))+O_p(\frac{h_2^2}{h})\\
&=O_p(h_2+\frac{h_2^{2}}{h}).
\end{split}
\end{eqnarray*}

In order to apply the theory for non-degenerate $U-$statistics\citep{serfling1980}, we also need to verified that $E\left(\|(B_{ij}+B_{ji})L\left(\frac{X_j-X_i}{h_0}\right)\|^2\right)=o(n).$
For any $p\times q$ matrix $A=(a_{ij})$, $\|A\|$ denote as its Euclidean norm, i.e., $\|A\|=[tr(AA^{\top})]^{1/2}.$

\begin{eqnarray*}
\begin{split}
&E\left(\|\frac{1}{2}(B_{ij}+B_{ji})\frac{1}{h_2}H\left(\frac{\alpha^{\top}X_j
-\alpha^{\top}X_i}{h_2}\right)\|^2\right)\\
&\leq 2E\left(\|\frac{1}{2}B_{ij}
\frac{1}{h_2^q}H\left(\frac{\alpha^{\top}X_j-\alpha^{\top}X_i}
{h_2}\right)\|^2\right)
+2E\left(\|\frac{1}{2}B_{ji}
\frac{1}{h_2^q}H\left(\frac{\alpha^{\top}X_j-\alpha^{\top}X_i}
{h_2^q}\right)\|^2\right)\\
&=E\left(\|B_{ij}\frac{1}{h_2^q}H\left(\frac{\alpha^{\top}X_j-\alpha^{\top}X_i}{h_0}\right)\|^2\right)\\
&=\frac{1}
{h_2^{2q}}E\left((s_p^1(X_{1i},\alpha^{\top}X_i))^2
K_i^2\frac{(D_j-p(X_i))^2}
{f^2(\alpha^{\top}X_i)}H^2
\left(\frac{\alpha^{\top}X_j-\alpha^{\top}X_i}{h_2}\right)\right)\\
&\qquad \qquad+\frac{1}
{h_2^{2q}}E\left((S_p^{(2)}(X_{1i},\alpha^{\top}X_i))^2K_i^2
\frac{(D_j-p(X_i))^2}
{f^2(X_i)}H^2\left(\frac{\alpha^{\top}X_j-\alpha^{\top}X_i}{h_2}\right)\right)
\end{split}
\end{eqnarray*}
where \begin{eqnarray*}
s_p(X_{1i},\alpha^{\top}X_i)&=&E(\phi_p(p(X_i),Z_i)\mid X_{1i},\alpha^{\top}X_i)
=\left(
\begin{array}{c}
 -\frac{E(\mathbb I(Y_i(1)^*\leq0)-\tau\mid X_{1i},\alpha^{\top}X_i)}{p(X_i)} \\
-\frac{E(\mathbb I(Y_i(1)^*\leq0)-\tau\mid X_{1i},\alpha^{\top}X_i)}{p(X_i)}\frac{X_{1i}-x_{10}}{h}
 \end{array}
\right)\\
&=&\left(
           \begin{array}{c}
            -\frac{\tilde m_{1,\tau}(X_{1i},\alpha^{\top}X_i)}{p(X_i)}\\
             -\frac{\tilde m_{1,\tau}(X_{1i},\alpha^{\top}X_i)}{p(X_i)}\frac{X_{1i}-x_{10}}{h}\\
           \end{array}
         \right)
:=\left(
           \begin{array}{c}
             s_p^1(X_{1i},\alpha^{\top}X_i) \\
             s_p^{(2)}(X_{1i},\alpha^{\top}X_i)\\
           \end{array}
         \right).
\end{eqnarray*}

It follows that
\begin{eqnarray}
\begin{split}
&\frac{1}
{h_2^{2q}}E\left((s_p^{(2)}(X_{1i},\alpha^{\top}X_i))^2K_i^2\frac{(D_j-p(\alpha^{\top}X_j))^2}
{f^2(\alpha^{\top}X_i)}H^2\left(\frac{\alpha^{\top}X_j-\alpha^{\top}X_i}{h_2}\right)\right)\\
&=\frac{1}
{h_2^{2q}}E\left((s_p^{(2)}(X_{1i},\alpha^{\top}X_i)))^2K_i^2\frac{p(\alpha^{\top}X_j)(1-p(\alpha^{\top}X_j))}
{f^2(\alpha^{\top}X_i)}H^2\left(\frac{\alpha^{\top}X_j-\alpha^{\top}X_i}{h_2}\right)\right)\\
&=\frac{1}
{h_2^q}\int s_p^{(2)}(X_{1i},\alpha^{\top}X_i)))^2K_i^2
\frac{p(\alpha^{\top}X_i+h_2u)(1-p(\alpha^{\top}X_i+h_2u))}
{f^2(\alpha^{\top}X_i)}H^2(u)dudX_i\\
&=O_p\left(\frac{1}{h_2^q}\right).
\end{split}
\end{eqnarray}
Under the assumption: $nh_2^q\rightarrow0$, we have  $\frac{1}
{h_2^{2q}}E\bigg((s_p^{(2)}(X_{1i},\alpha^{\top}X_i))^2
K_i^2\frac{(D_j-p(\alpha^{\top}X_j))^2}
{f^2(\alpha^{\top}X_i)}\\
\times H^2\left(\frac{\alpha^{\top}X_j-\alpha^{\top}X_i}{h_2}\right)\bigg)
=o(n).$ Similarly we can verify that $\frac{1}
{h_2^{2q}}E\bigg((s_p^{(1)}(X_{1i},\alpha^{\top}X_i))^2\\
\times K_i^2\frac{(D_j-p(\alpha^{\top}X_j))^2}
{f^2(\alpha^{\top}X_i)}H^2\left(\frac{\alpha^{\top}X_j-\alpha^{\top}X_i}{h_2}\right)\bigg)
=O_p\left(\frac{1}{h_2^q}\right)=o(n).$
Thus the condition of Lemma 3.1 of \cite{zheng1996} is satisfied. Further, combining the fact $\varepsilon_j$' s are
mutually independent conditional on the sample path of the
$X_j$' s, we can obtain the following result, with  assumptions: $h_2\rightarrow 0$ and ${h_2^{2s_2}}{h^{-(2s_2+1)}}\rightarrow 0$ and $ s_2 \geq 2$,
\begin{eqnarray}\label{j11}
& &\tilde J_{11}=\frac{1}{\sqrt{nh}}\sum_{i=1}^nO_p(h_2+\frac{h_2^{2}}{h})\epsilon_i+o_p(1)\nonumber\\
& &\qquad=\frac{1}{\sqrt{n}}\sum_{i=1}^nO_p(\frac {h_2}{\sqrt h}+\frac{h_2^{2}}{h^{3/2}})\epsilon_i+o_p(1)=o_p(1),
\end{eqnarray}
which is similar as the proof for $J_{12}$ for Theorem 2.3. in the nonparametric setting.

{ Combined with  $\tilde J_{13}=o_p(1)$, $\tilde J_{2}=o_p(1)$ and $\tilde J_{3}=o_p(1)$}, it yields
the asymptotically linear representation of $\hat q^{scqte}_{1,\tau}(x_{10})$  with the semiparametric estimator  $\hat p(\alpha^{\top}X)$ is
\begin{eqnarray*}
\sqrt{nh}\left(\hat q^{scqte}_{1,\tau }(x_{10})- {q_{1,\tau}
(x_{10})}\right)=-\frac{1}{\sqrt{nh}}
\frac{1}{f(x_{10})}\sum_{i=1}^{n}
\left(\frac{D_i}{p(\alpha^{\top}X_i)} \eta_{1,
\tau}(Y_i)\right)K_i+o_p(1),
\end{eqnarray*}

Similarly, the asymptotical linear representation of $\hat
q^{scqte}_{0,\tau}(x_{10})$ is as
\begin{eqnarray*}
\sqrt{nh}\left(\hat q^{scqte}_{0,\tau }(x_{10})- {q_{0,\tau}
(x_{10})}\right)=-\frac{1}{\sqrt{nh}}
\frac{1}{f(x_{10})}\sum_{i=1}^{n}
\left(\frac{D_i}{1-p(\alpha^{\top}X_i)} \eta_{0,
\tau}(Y_i)\right)K_i+o_p(1),
\end{eqnarray*}
where $K_i=K(\frac{X_{1i}-x_{10}}{h})$,
$\eta_{j,\tau}(Y_i)=\frac{\mathbb I(Y_i\leq \tilde
q_{j,\tau}(X_{1i}))-\tau}
{f_j(q_{j,\tau}(x_{10})\mid x_{10})}$, $j=0,1.$

Thus we can get the asymptotical linear representation of
$\hat\Delta^{scqte}_\tau(x_{10})$ as
\begin{eqnarray*}
& &\sqrt{nh}\left(\hat \Delta^{scqte}_{\tau }(x_{10})-
{\Delta_{\tau} (x_{10})}\right)=-\frac{1}{\sqrt{nh}}
\frac{1}{f(x_{10})}\sum_{i=1}^{n}\phi_3(p(X_i),Z_i)K_i+o_p(1),
\end{eqnarray*}
where $\phi_3(p(X_i),Z_i)=\phi_1
(p(X_i),Z_i)$. Further, $\phi_1
(p(X_i),Z_i)=
\frac{D_i}{p(\alpha^{\top}X_i)} \eta_{1,
\tau}(Y_i) -\frac{1-D_i}{1-p(\alpha^{\top}X_i)} \eta_{0,\tau}(Y_i)$,
$\eta_{j,\tau}(Y_i)=\frac{\mathbb I(Y_i\leq
\tilde q_{j,\tau}(X_{1i}))-\tau}
{f_j(q_{j,\tau}(x_{10})\mid x_{10})},j=0,1.$

{ More generally, when $1<r\leq l$ and $X_1 \sqsubset^{l-r}\alpha^{\top}X$, we can similarly prove
\begin{eqnarray}\label{j11}
& &\tilde J_{11}=\frac{1}{\sqrt{nh^l}}\sum_{i=1}^nO_p(h_2^{r}+\frac{h_2^{r+1}}{h})\epsilon_i+o_p(1)\nonumber\\
& &\qquad=\frac{1}{\sqrt{n}}\sum_{i=1}^nO_p(\frac {h_2^{r}}{  h{^{l/2}}}+\frac{h_2^{r+1}}{h^{1+l/2}})\epsilon_i+o_p(1).
\end{eqnarray}
 Under assumptions $9'$ with $s_2(2-l/r)>0$, it follows that $\tilde J_{11}=o_p(1)$.
 Thus, by the analogy of the discussion of $X_1\in R$ and $X_1\sqsubset^0\alpha^{\top}X$, we complete the proof of Lemma \ref{out_lemma}.
}

 \hfill$\Box$

\subsubsection{Proof of {Lemma~S.2}.}

Note that if $X_1 \subseteq \alpha^{\top}X$, $s_p(X_1,\alpha^{\top}X)=s_p(\alpha^{\top}X)$. It follows that
\begin{eqnarray*}
& &\frac{1}{f(\alpha^{\top}X_i)}E\bigg(\frac{1}{h_2^q}H(\frac{\alpha^{\top}X_j-\alpha^{\top}X_i}{h_2})
s_p(\alpha^{\top}X_j)K(\frac{X_{1j}-x_{10}}{h})\mid \alpha^{\top}X_i\bigg)\\
& &=s_p(\alpha^{\top}X_i)K(\frac{X_{1i}-x_{10}}{h})+O_p(\frac{h_2^{ s_2}}{h^{s_2}}).
\end{eqnarray*}
where $s_p(\alpha^{\top}X)=-\frac{E(\mathbb I(y(1)\leq \tilde q_{1,\tau}(X_{1}))-\tau\mid \alpha^{\top}X)}{p(\alpha^{\top}X)}.$

Thus, the proof of $Lemmay$ S.2 can be similar to that of $NCQTE(x_{10})$ in Theorem~2.3 to  provide a smaller variance.
 \hfill$\Box$

\subsection{Proof of Corollary~2.2.}

 From the proof for Theorem~2.3, we can see that if we can further derive that if the term $J_1$ in (\ref{j1}) goes to zero in probability, $NCQTE(x_{10})$ will have the same limiting distribution as $PCQTE(x_{10})$. By the condition $\sqrt {nh}(h_0^{s}+\sqrt {\log (n)/nh_0^k})\to 0$, we can derive that
$$
|J_1|\le O_p(\sqrt {nh}(h_0^{s}+\sqrt {\frac {\log (n)}{nh_0^k}}))\frac 1{nh}\sum_{j=1}^n| S_p(X_j)
K_j|=o_p(1),
$$
because $\frac 1{nh}\sum_{j=1}^n| S_p(X_j)
K_j|=O_p(1).$

For $SCQTE(x_{10})$, we can similarly deal with the term $\tilde J_1$ in (\ref{tj1}) such that the limiting distribution can be the same as that of $PCQTE(x_{10})$.  The proof is finished.
 \hfill$\Box$

\subsubsection{Proof of Corollary~2.3.}
 By the condition, $D\perp X\mid \tilde X$ with $\tilde X\subseteq X$. we can immediately derive that the propensity score $p(X)=p(\tilde X)$   as $p(X)$ is the conditional probability of $D=1$ when $X$ is given. From the proof for Theorem~2.3,  we only need to prove that the term  $J_1=o_p(1)$. Let $\tilde X\in R^{d}$ with $d< k$. Without loss of generality, assume $d=k-1$ and  $X=(X_1,\tilde X)$ and $\hat p(\tilde X_i)=\sum_{j\neq i}w_{ij}D_j,$ where
$w_{ij}=\frac{L\left(\frac{\tilde X_j-\tilde X_i}{h_0}\right)}{\sum_{j\neq
i}L\left(\frac{\tilde X_j-\tilde X_i}{h_0}\right)}.$
By the analogy discussion of $Lemma~s.1$, we consider the case $l=1$, i.e. $X_1\sqsubset^{0}\tilde X$ at first.

We can decompose $\hat
p(X_i)$ as
\begin{eqnarray}
\hat p(\tilde X_i)-p(\tilde X_i)=\sum_{j:j\neq
i}w_{ij}\epsilon_j+\vartheta_n(X_i),
\end{eqnarray}
where $\epsilon_j=D_j-p(\tilde X_j),$ $\vartheta_n(\tilde X_i)=\sum_{j\neq
i}w_{ij}p(\tilde X_j)-p(\tilde X_i).$
 The corresponding properties of $\hat p(\tilde X)$ is
\begin{eqnarray}
&&\sup_{\tilde X \in \tilde \chi} |\hat p(\tilde X)-p(\tilde X)|=O_p\left(h_0^s+\sqrt{\frac{\log n}{nh_0^d}}\right),\nonumber\\
&&\sup_{\tilde X \in \tilde \chi} |\vartheta_n(\tilde X)|=O_p\left(h_0^s\right).
\end{eqnarray}
Similarly as the proof for Theorem~2.4, $J_1$ can be decomposed as
\begin{eqnarray}\label{j1}
& &J_1=\frac{1}{\sqrt{nh}}\sum_{i=1}^{n}S_p(X_i) K_i(\hat
p(\tilde X_i)- p(\tilde X_i))\nonumber\\
& &=\frac{1}{\sqrt{nh}}\sum_{i=1}^n\sum_{j:j\neq i}\epsilon_jw_{ij}S_p(X_i)K_i
+\frac{1}{\sqrt{nh}}\sum_{i=1}^{n}S_p(X_i)
K_i\vartheta_n(\tilde X_i)\nonumber\\
& &:=\overline J_{11}+J_{13}.
\end{eqnarray}
Here $S_p(X_i)=E(\phi_p(p(\tilde X_i),Z_i)\mid X_i).$
The similar treatment for $\tilde J_{11}$ in $SCQTE(x_{10})$, under the assumptions $h_0\rightarrow 0$,  $nh_0^d\rightarrow \infty$, we can also rewrite $\overline J_{11}$ as
\begin{eqnarray*}
\overline J_{11}=\frac{n-1}{\sqrt{nh}}\frac{1}{n(n-1)}
\sum_{i=1}^n\sum_{j\neq i}^n[\frac{H_{ij}+H_{ji}}{2}]\frac{1}{h_0^d}
L(\frac{\tilde X_j-\tilde X_i}{h_0})+o_p(1).
\end{eqnarray*}
Here $
H_{ij}=S_p(X_i) K_i\frac{1}
{h_0^df(\tilde X_i)}\epsilon_j.$
We can also show that $E\left(\|(H_{ij}+H_{ji})\frac{1}{h_0^d}L\left(\frac{X_j-X_i}{h_0}\right)\|^2\right)
=O_p\left(\frac{1}{h_0^d}\right)=o(n)$,
 and $E(H_{ij}\mid Z_i)=0$ and
 \begin{eqnarray*}
 \begin{split}
 &E\left(H_{ji}\frac{1}{h_0^d}L\left(\frac{\tilde X_j-\tilde X_i}{h_0}\right)\mid Z_i\right)\\
&=E\left(S_p(X_{j} )K(\frac{x_{1j}-x_{10}}{h})
\frac{1}{f(\tilde X_j)}\frac{1}{h_0^d}L(\frac{\tilde X_j-\tilde X_i}{h_0})\mid X_i\right)\varepsilon_i.\\
\end{split}
\end{eqnarray*}
Let $X_{1}=\nu_1,\tilde X=\nu_2$, and denote $(\frac{\nu_1-\nu_{1i}}{h_0},\frac{\nu_2-\nu_{2i}}{h_0})$ as $(t_1,t_2)$. Similarly as the proof in $SCQTE(x_{10})$, we have
\begin{eqnarray*}
&&E\left(S_p(X_{j} )K(\frac{x_{1j}-x_{10}}{h})
\frac{1}{f(\tilde X_j)}\frac{1}{h_0^d}L(\frac{\tilde X_j-\tilde X_i}{h_0})\mid X_i\right)\nonumber\\
&=&h_0\int K(\frac{t_1h_0}{h}+
\frac{\nu_{1i}-x_{10}}{h})\frac{1}
{f(\nu_{2i}+h_2t_2)}
L(t_2)S_p(\nu_{1i}+h_0t_1,\nu_{2i}+h_0t_2)\\
&&
f_{\nu_{12}}(\nu_{1i}+h_0t_1,\nu_{2i}+h_0t_2)
dt_1dt_2\nonumber\\
&=&O_p(h_0)+O_p(\frac{h_0^2}{h}),
\end{eqnarray*}
where $f_{\nu_{12}}(\nu_1,\nu_2)$ is the joint density function of $(X_{1},\tilde X).$
Similarly as that in (\ref{j11}),
under the  assumptions: $h_0\rightarrow 0$ and $\frac{h_0^{2s}}{h^{2s+1}}\rightarrow 0$, we have $\overline J_{11}=o_p(1).$ The other terms  $J_{13}$, $J_{2}$ and $J_3$ can be similarly proved to be $o_p(1).$ While $1<r\leq l$, under the additional condition $s(2-l/r)>0$, we can also prove that $\overline J_{11}=o_p(1)$.
 Thus we complete the proof of this corollary.
 \hfill$\Box$

\subsubsection{Proof of Corollary~2.4.}

Decompose $m_{j,\tau}(X)=m_{j,\tau}(\alpha^{\top}X)+\Delta m_{j,\tau}$, $j=0,1$. Via using the fact that $p(X)=p(\tilde X)=p(\alpha^{\top}X)$ and some direct calculations, we can get
\begin{eqnarray*}
 \sigma_{scqte}^{*2}(x_{10})=\sigma_{ncqte}^{*2}(x_{10})+
E\bigg[p(\alpha^{\top}X)(1-p(\alpha^{\top}X))\big(\frac{\Delta m_{1,\tau}}{p(\alpha^{\top}X)}
 +\frac{\Delta m_{0,\tau}}{1-p(\alpha^{\top}X)}\big)^2\mid X_1=x_{10}\bigg]
\geq \sigma_{ncqte}^{*2}(x_{10}).
\end{eqnarray*}

 \hfill$\Box$

\bibhang=1.7pc
\bibsep=2pt
\fontsize{9}{14pt plus.8pt minus .6pt}\selectfont
\renewcommand\bibname{\large \bf References}
\expandafter\ifx\csname
natexlab\endcsname\relax\def\natexlab#1{#1}\fi
\expandafter\ifx\csname url\endcsname\relax
  \def\url#1{\texttt{#1}}\fi
\expandafter\ifx\csname urlprefix\endcsname\relax\def\urlprefix{URL}\fi

\lhead[\footnotesize\thepage\fancyplain{}\leftmark]{}\rhead[]{\fancyplain{}\rightmark\footnotesize{} }